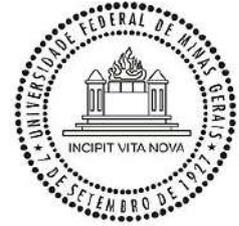

*Célio Augusto Terra de Souza*

# Dynamic Phenomena in Interacting Particle Systems: Phase Transitions and Equilibrium

Célio Augusto Terra de Souza

# Dynamic Phenomena in Interacting Particle Systems: Phase Transitions and Equilibrium

Thesis submitted to the Graduate Program of Mathematics of the Universidade Federal de Minas Gerais in partial fulfillment of the requirements for the degree of *Philosophiæ Doctor* in Mathematics.

Advisor: Prof. Dr. Leonardo Trivellato Rolla

Co-advisor: Prof. Dr. Bernardo Nunes Borges de Lima

Belo Horizonte - MG

2024

ATA DA DUCENTÉSIMA VIGÉSIMA TERCEIRA DEFESA DE TESE DE DOUTORADO DO ALUNO CÉLIO AUGUSTO TERRA DE SOUZA, REGULARMENTE MATRICULADO NO PROGRAMA DE PÓS-GRADUAÇÃO EM MATEMÁTICA DO INSTITUTO DE CIÊNCIAS EXATAS DA UNIVERSIDADE FEDERAL DE MINAS GERAIS, REALIZADA DIA 21 DE NOVEMBRO DE 2024.

Aos vinte e um dias do mês de novembro de 2024, às 10h00, na sala 3060, reuniram-se os professores abaixo relacionados, formando a Comissão Examinadora homologada pelo Colegiado do Programa de Pós-Graduação em Matemática, para julgar a defesa de tese do aluno **Célio Augusto Terra de Souza**, intitulada: "*Dynamic Phenomena in Interacting Particle Systems: Phase Transition and Equilibrium*", requisito final para obtenção do Grau de Doutor em Matemática. Abrindo a sessão, o Senhor Presidente da Comissão, Prof. Leonardo Trivellato Rolla, após dar conhecimento aos presentes do teor das normas regulamentares do trabalho final, passou a palavra ao aluno para apresentação de seu trabalho. Seguiu-se a arguição pelos examinadores com a respectiva defesa do aluno. Após a defesa, os membros da banca examinadora reuniram-se reservadamente, sem a presença do aluno, para julgamento e expedição do resultado final. Foi atribuída a seguinte indicação: o aluno foi considerado aprovado sem ressalvas e por unanimidade. O resultado final foi comunicado publicamente ao aluno pelo Senhor Presidente da Comissão. Nada mais havendo a tratar, o Presidente encerrou a reunião e lavrou a presente Ata, que será assinada por todos os membros participantes da banca examinadora. Belo Horizonte, 21 de novembro de 2024.

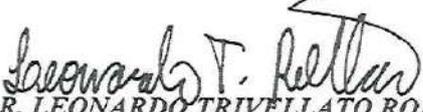

**PROF. DR. LEONARDO TRIVELLATO ROLLA**
Orientador (USP)

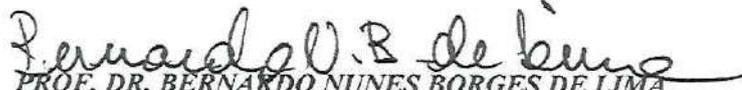

**PROF. DR. BERNARDO NUNES BORGES DE LIMA**
Coorientador (UFMG)

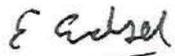

**PROF. DR. ENRIQUE DANIEL ANDJEL**
Examinador (Universidade de Marselha)

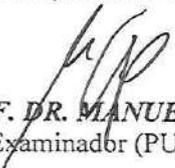

**PROF. DR. MANUEL CABEZAS**
Examinador (PUC-Chile)

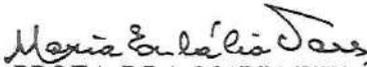

**PROFA. DRA. MARIA EULÁLIA VARES**
Examinadora (UFRJ)

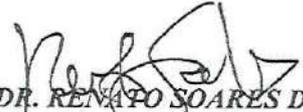

**PROF. DR. RENATO SOARES DOS SANTOS**
Examinador (UFMG)

# Acknowledgements

The journey to this moment has spanned four years of hard work, and I could not have done it alone. I would like to express my heartfelt gratitude to those who supported me along the way.

First and foremost, I thank God, the foundation of our lives.

I am deeply grateful to my family—my father, mother, brother, grandparents, and everyone who has been and will always be my safe haven.

I extend my appreciation to Léo and Bernardo, who believed in me and fought tirelessly to help me learn and grow. Thank you for not giving up even when you had to repeat something many times until I finally understood.

My sincere thanks go to Professor Pablo Groisman for our collaborative work on part of this thesis and for his invaluable feedback.

I would also like to thank Andréa and Kelli for their assistance with the myriad of bureaucratic tasks and for the moments of relaxation and camaraderie in the break room.

Last but certainly not least, I am thankful for my friends. You know how much you mean to me and how vital your support has been. I won't name names, but you all know who you are. I could never have made it this far without you. Thank you so much.

This study was financed in part by the Coordenação de Aperfeiçoamento de Pessoal de Nível Superior – Brasil (CAPES) – Finance Code 001, and by Fundação de Amparo à Pesquisa do Estado de Minas Gerais (FAPEMIG).

*This thesis is dedicated*
*to everyone who contributed*
*in these four years*
*and made this work possible.*

*The thoughts of pure mathematics are true, not approximate or doubtful;*
*they may not be the most interesting or important of God's thoughts,*
*but they are the only ones that we know exactly.*
Hilda Hudson

*"Well, in our country," said Alice, still panting a little, "you'd generally get*
*to somewhere else—if you ran very fast for a long time, as we've been doing."*
*"A slow sort of country!" said the Queen. "Now, here, you see, it takes all*
*the running you can do, to keep in the same place. If you want to get*
*somewhere else, you must run at least twice as fast as that!"*
Lewis Carroll, *Through the looking-glass and what Alice found there*

# Resumo

Esta tese investiga fenômenos críticos e estados de equilíbrio em diversos modelos estocásticos por meio de três estudos interligados.

No primeiro capítulo, analisamos o modelo de Passeios Aleatórios Ativados em um anel unidimensional no regime de alta densidade. Introduzimos um procedimento de *topplings* que constrói incrementalmente um ambiente que demonstra a atividade sustentada por longos períodos. Esta abordagem fornece uma prova concisa e auto-contida da existência de uma fase lenta para taxas de sono arbitrariamente grandes.

O segundo capítulo concentra-se em um processo de contato unidimensional modificado com taxas de infecção distintas. Especificamente, a infecção se espalha a uma taxa $\lambda_e$ nos limites da região infectada e a uma taxa $\lambda_i$ em outros lugares. Estabelecemos a existência de uma medida invariante para este processo quando $\lambda_i = \lambda_c$, $\lambda_e = \lambda_c + \varepsilon$ onde $\lambda_c$ denota o parâmetro crítico para o processo de contato padrão. Além disso, demonstramos que o processo, quando observado pela borda direita, converge fracamente para esta medida invariante. Mostramos também que a infecção morre quase certamente ao longo da curva crítica dentro da região atrativa do espaço de fase.

No capítulo final, exploramos distribuições quase-estacionárias (DQE) para dois processos populacionais subcríticos em tempo contínuo: passeios aleatórios com ramificação e processos de ramificação com genealogia. Provamos a existência e a unicidade da DQE para esses processos, aproveitando os aspectos espaciais de sua dinâmica.

**Palavras-chave:** fenômenos críticos; passeios aleatórios ativados; processo de contato com fronteira modificada; convergência para equilíbrio; processos populacionais; distribuições quase-estacionárias.

# Abstract

This thesis investigates critical phenomena and equilibrium states in various stochastic models through three interconnected studies.

In the first chapter, we analyze the Activated Random Walk model on a one-dimensional ring in the high-density regime. We introduce a toppling procedure that incrementally constructs an environment demonstrating the sustained activity over extended periods. This approach provides a concise and self-contained proof of the existence of a slow phase for arbitrarily large sleep rates.

The second chapter focuses on a modified unidimensional contact process with varying infection rates. Specifically, infection spreads at rate $\lambda_e$ at the boundaries of the infected region and at rate $\lambda_i$ elsewhere. We establish the existence of an invariant measure for this process when $\lambda_i = \lambda_c$, $\lambda_e = \lambda_c + \varepsilon$ where $\lambda_c$ denotes the critical parameter for the standard contact process. Furthermore, we demonstrate that the process, when observed from the right edge, converges weakly to this invariant measure. We also show that infection dies almost surely along the critical curve within the attractive region of the phase space.

In the final chapter, we explore quasi-stationary distributions (QSDs) for two subcritical population processes in continuous time: branching random walks and branching processes with genealogy. We prove the existence and uniqueness of QSDs for these processes by leveraging spatial aspects of their dynamics.

**Keywords:** critical phenomena; activated random walks; modified boundary contact process; convergence to equilibrium; populational processes; quasi-stationary distributions.

# List of Figures

# Contents

# Introduction

The study of interacting particle systems is a cornerstone of modern statistical mechanics and probability theory, offering profound insights into the behavior of complex systems. These systems, characterized by the interactions among numerous individual components, exhibit emergent behaviors that can often be counterintuitive. Over the past several decades, a diverse array of models has been developed to explore different facets of these systems, leading to a rich and interdisciplinary body of research. Prominent models, such as the voter model, the Ising model, the contact process, and the frog model, have been pivotal in advancing our understanding of a range of phenomena, from social dynamics and opinion formation to physical phase transitions and critical phenomena.

Interacting particle systems have a wide range of applications across various fields, reflecting their fundamental role in understanding complex phenomena. In statistical physics, they model critical phenomena and phase transitions, such as in the Ising model for ferromagnetism or the contact process for disease spread. Biology benefits from these models to study population dynamics, including the spread of infectious diseases and the dynamics of ecosystems, where models like the frog model or branching processes provide insights into survival and extinction patterns. In epidemiology, interacting particle systems help simulate the spread of diseases, allowing researchers to explore how infections propagate through populations and evaluate the impact of public health interventions. Social sciences use these models to understand opinion formation, rumor spreading, and network dynamics, providing a framework to analyze how individuals influence each other in social networks. Additionally, in computer science, interacting particle systems contribute to algorithms for distributed computing and optimization problems, where they model resource allocation and system performance. These diverse applications illustrate the versatility and importance of interacting particle systems in both theoretical and practical contexts.

This thesis aims to address some questions related to phase transitions and convergence to equilibrium in interacting particle systems. By examining these systems under varying conditions and constraints, we seek to deepen our understanding of their underlying principles and behaviors.

# Phase Transition in Activated Random Walks

Activated Random Walks (ARWs) represent a fascinating extension of traditional random walk models by incorporating activation-inactivation mechanisms. In the ARW framework, there are two types of particles in a lattice. *Active* particles perform independent continuous-time random walks at rate 1 and each active particle becomes *sleeping* at rate $\lambda$, and remain sleeping until an active particle reaches the site where the sleeping particle is. Sleeping particles do not move. This added layer of complexity introduces a wealth of dynamical behavior and allows for the exploration of novel phenomena, such as phase transitions and criticality, within the context of random walks.

A particularly intriguing aspect of the ARW model is its connection to the concept of *self-organized criticality*. Self-organized criticality, introduced in [BTW87], describes systems that naturally evolve into a critical state without requiring fine-tuning of external parameters. In these systems, criticality arises spontaneously due to the interactions among the system's components. Besides Activated Random Walks, other models that exhibits self-organized criticality are the Manna Sandpile Model [Man91] and its Abelian variant, the Stochastic Sandpile Model [Dha90].

The existence of self-organized criticality is closely related to the occurrence of *phase transitions* in models with explicit parameters. For ARWs, this phase transition manifests in the *fixed-energy* model. In this conservative model, particles are neither created nor destroyed, and the key parameter of interest is the density of particles. The fixed-energy model can be studied on both infinite lattices and finite torii. In Chapter 1, we provide a self-contained proof of the existence of a phase transition in the fixed-energy model on a one-dimensional ring. In this case, the system will reach an absorbing configuration if, and only if the number of particles is less than or equal to $N$. We prove that in the ARW model in the onedimensional torus $\mathbb{Z}_N$, for every sleeping rate the time for absorption is exponential in $N$ if the number of particles is close enough to $N$ (but smaller than $N$) (Theorem 1.1). These results, which have been published in [LRT24], contribute to the broader understanding of phase transitions in stochastic systems.

# Convergence to Equilibrium in the Modified Boundary Contact Process

The contact process is a foundational model in the study of interacting systems, shedding light on essential behaviors in population dynamics, such as the spread of infections, rumors, and social behaviors. Originally proposed by Harris in 1974 [Har74], this model features a population where individuals can be in one of two states: *infected* or *healthy*. Infected individuals can transmit the infection to their healthy neighbors, while they recover at a certain rate. This interaction framework mirrors the complexity of real-world epidemic dynamics and social contagion phenomena. Understanding how such systems evolve towards stable distributions, or invariant measures, is crucial for grasping their long-term behavior.

Several modifications to the classical contact process have been introduced to explore different aspects of its dynamics. One notable variation, discussed in [DS00], involves distinct infection rates for the boundary and the interior of the infected region. The phase space of this boundary-modified model was initially explored in [DS00] and further examined in [AR23], which also established convergence of the process seen from the edge to an equilibrium state within the attractive region of the phase space. In Chapter 2, we build on these findings by extending the results of convergence to equilibrium to the non-attractive region (Theorem 2.2) and refining the phase space characterization. Specifically, we demonstrate that infection almost surely dies out in the attractive segment of the critical curve (Theorem 2.5).

# Quasi-Equilibrium in Population Processes

Quasi-stationary distributions (QSDs) are fundamental to understanding population processes, offering deep insights into the long-term dynamics of systems that face eventual extinction. A QSD is a distribution that is invariant for the process conditioned that extinction has not yet occurred. QSDs are particularly significant in scenarios where populations experience random fluctuations and face the possibility of extinction, such as in biological models, epidemiological studies, and ecological systems.

In various subcritical population models, whether discrete-time or continuous-time Markov chains, the concept of a quasi-stationary distribution provides a refined understanding of the system's behavior before extinction. Unlike stationary distributions, which describe the long-term behavior of a process that continues indefinitely, QSDs focus on the behavior conditioned on survival. They offer insights into how the population is distributed across different states given that it has not yet faced extinction. For example, QSDs play a role in describing scaling limits for subcritical contact processes [DR17] and can also represent the limit of the distribution observed from the edge of a process [AEGR15].

The behavior of quasi-stationary distributions (QSDs) can differ significantly across various models. For example, in subcritical branching processes, multiple QSDs exist, each representing a combination of the minimal QSD (the one with the highest absorption rate) and a QSD specific to the pure-death process [Cav78]. In contrast, subcritical contact processes possess a unique QSD [AGR20]. This difference may be attributed to the geometric aspects inherent in the contact process that are absent in branching processes. In Chapter 3, we investigate two variants of classical branching processes that incorporate geometric elements: the branching process with genalogy, that gives information about the genealogical relationships of the individuals, and the branching random walks, that adds information about the spatial location of the individuals. We then prove existence and uniqueness of quasi-stationary distribution for those two processes (Theorems 3.3 and 3.4). This chapter provides new insights into the intricate interplay between geometry and quasi-stationarity in population processes. This chapter was done in collaboration with Pablo Groisman.

# Part I

# Phase transition in finite systems

# Chapter 1

# Slow phase for the fixed energy ARW

Activated Random Walk (ARW) is an Abelian model of interacting particles that can be described as follows. Particles are distributed on a transitive graph. There are two types of particles, *active* and *sleeping* particles. Active particles perform independent continuous-time random walks at rate 1, and fall asleep at rate $\lambda$. Sleeping particles do not move, and continue to sleep until an active particle reaches the same site, and then become active again. See [Rol20] for an introduction and main results. Recent progress in the model has been made, including bounds for the critical density in the infinite lattice [AFG24, ARS22, FG22, Hu22, Tag23], and of related models [HHRR22, PR21], separation cutoff for the total variation distance to equilibrium [BS22], bounds on the mixing time of the driven-dissipative model [LL21], scaling limits [CR21], necessary conditions for supercriticality [For24], bounds for the spread of the particles [LS21] and a prove of the universality conjecture [HJJ24].

In this work, we will focus on ARWs on the unidimensional ring $\mathbb{Z}_N := \mathbb{Z}/N\mathbb{Z}$. We suppose a deterministic initial configuration with $\zeta N$ particles. As the ring $\mathbb{Z}_N$ is finite, it is immediate that the system will eventually be absorbed in a configuration with only sleeping particles (i.e., will *stabilize*) if, and only if, $\zeta \leq 1$. Thus, we are interested in *how much time* time the system will take to stabilize. Our main result is the following theorem, which has been published in [LRT24].

**Theorem 1.1.** *Let $\mathcal{J}$ be the total number of jumps the particles do until the*

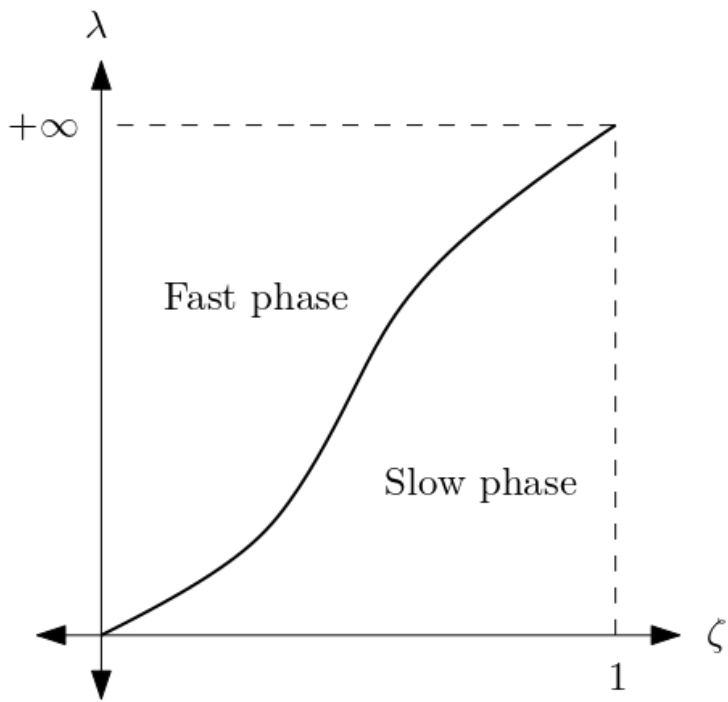

**Figure 1.1:** Phase diagram for the unidimensional ring

*initial configuration is stabilized. For every $0 < \lambda < +\infty$, there are constants $\delta, \delta' > 0$, depending on $\lambda$, such that, if $\zeta$ is close enough to 1,*

$$\mathbb{P}(\mathcal{J} \geq e^{\delta N}) \geq 1 - e^{-\delta' N},$$

*for every $N$ sufficiently large.*

As the existence of a polynomial phase (i.e., a fast phase) for small densities was already proven in [RS12, BGHR19], the above result shows that there is a phase transition on the ring $\mathbb{Z}_N$. Thus, we have the phase diagram as showed in Figure 1.1.

While this result was being prepared, two papers came out. In [AFG24], Theorem 1.1 is proved for the $d$-dimensional torus, $d \geq 2$, using auxiliary results and arguments of [FG22]. Their proof uses a hierarchical argument. Fixed a "dormitory" set where the particles will eventually stabilize, the model is reduced to a model of density 1 in the dormitory, which is then divided in a series of growing "clusters". Stabilization of each of those clusters will affect the configuration in the other clusters, leading to a long stabilization time with high probability, even after summing over all the possible dormitory sets. A modification of this argument shows that the result is still true in the one-dimensional case. In [HJJ24] it is shown that the critical density on $\mathbb{Z}_N$ and in an infinite version of ARW in $\mathbb{Z}$ are the same. Their proof is based on the cyclical argument of [BGHR19], using a "percolation in layers" approach to show that, in every cycle, there will remain a fraction of active particles that will be sufficient to wake up the sleeping particles.

Our proof uses a different approach which is self-contained and shorter. The main strategy is to build a "carpet" where active particles can perform long excursions. Unlike [HRR23], we do not assume a convenient initial configuration, instead our carpet is built on the fly while the system dynamics develop. We then divide the ring in large blocks, and prove that with high probability the system reaches certain configurations having many active particles in a certain region of the ring. In contrast to [BGHR19], which handles the low $\lambda$ regime, at the end of each cycle we are left with inconvenient configurations connecting the source and sink regions, and need to make do with them.

This chapter is structured as follows. In Section 1.1 we give a formal description of the ARW model. In Section 1.2 we describe a dynamic procedure to move particles on the ring according to the ARW rules. This procedure will conserve a series of properties of the configuration, as proved in Section 1.3. Theorem 1.1 is proved in Section 1.4, using the cycle argument. In Sections 1.5 to 1.7 some auxiliaries results and lemmas are proved.

## 1.1 Description of the Model

In this section we give a more formal description of the model. The configuration space will be $\{0, \mathfrak{s}, 1, 2, \dots\}^{\mathbb{Z}_N}$, where $\mathfrak{s}$ is a symbol used to denote one sleeping particle. We define $\mathfrak{s} + 1 = 2$.

To each site $x$ on $\mathbb{Z}_N$ is associated a *stack of instructions* $\xi^x := (\xi^x_k)_{k \in \mathbb{N}}$ sampled independently over $x$ and $k$ accordingly to the following distribution

$$\xi^x_k = \begin{cases} \mathfrak{t}_{x,x+1} & \text{with probability } \frac{1}{2(1+\lambda)}; \\ \mathfrak{t}_{x,x-1} & \text{with probability } \frac{1}{2(1+\lambda)}; \\ \mathfrak{t}_{x,\mathfrak{s}} & \text{with probability } \frac{\lambda}{1+\lambda}, \end{cases}$$

where each one of the above transition acts on a configuration $\eta \in \{0, \mathfrak{s}, 1, 2, \dots\}^{\mathbb{Z}_N}$ in the following way

$$[\mathfrak{t}_{x,y}\eta](z) = \begin{cases} \eta(z) + 1 & \text{if } z = y; \\ \eta(z) - 1 & \text{if } z = x \\ \eta(z) & \text{otherwise,} \end{cases}$$

and

$$[\mathfrak{t}_{x,\mathfrak{s}}\eta](z) = \begin{cases} \mathfrak{s}, & \text{if } z = x \text{ and } \eta(x) = 1; \\ \eta(z), & \text{otherwise.} \end{cases}$$

We also introduce the *odometer field* $h = (h(x); x \in \mathbb{Z}_N)$, which count the number of instructions that have been utilized at each site site. The initial value of the odometer is $h \equiv 0$.

A site $x$ is called *stable* in a configuration $\eta$ if $\eta(x) = 0$ or $\mathfrak{s}$, and *unstable* otherwise. A configuration $\eta$ is *stable* if every $x$ is stable in $\eta$. The operation of *toppling* a site $x$, denoted by $\Phi_x$, is defined by

$$\Psi_x(\eta, h) = (\mathfrak{t}^x_{h(x)+1}, h + \delta_x).$$

In other words, to topple a site is to update the configuration by making a particle at site $x$ to obey the first instruction at site $x$ which have note been used. The toppling of a site $x$ is *legal* if $x$ is unstable.

A *legal sequence of topplings* is a $k$-tuple $\alpha = (x_1, x_2, \dots, x_k)$ of sites where $x_1$ is unstable, $x_2$ is unstable after $x_1$ is toppled, $x_3$ is unstable after $x_1$ and $x_2$ are toppled, etc. The *odometer of $x$ in $\alpha$*, denoted by $m_\alpha(x)$ is the number of times the site $x$ appears in $\alpha$. A legal sequence $\alpha$ *stabilizes* $\eta$ if after we topple all the sites of $\alpha$ we obtain a stable configuration.

The *Abelian property* of the ARW [Rol20, Lemma 2.4] states that if $\alpha$ and $\beta$ are two legal sequences of topplings that stabilize a configuration $\eta$, then $m_\alpha = m_\beta$. In other words, the order in which legal topplings are performed does not affect the final configuration obtained.

For a configuration $\eta$, we define the *odometer* of $\eta$ to be $m_\eta = m_\alpha$, with $\alpha$ a legal sequence of topplings that stabilize $\eta$. By the Abelian property, this definition does not depend on our choice of $\alpha$.

## 1.2 The toppling procedure

In this section, we will describe a sequence of topplings in the ring $\mathbb{Z}_N$. In Section 1.4, we will use this procedure as the main step of an iterative procedure to stabilize the configuration on the closed ring.

### 1.2.1 Preliminary steps

Without loss of generality, we can assume that the initial configuration has at most one particle per site. Otherwise, because $\zeta < 1$, we can a.s. topple sites with more than one particle until there are no sites with more than one particle.

Let $a \in 2\mathbb{N}$ be fixed, its value will be determined later, and let $K := a^2$. For simplicity, we assume $N = (n + 2)K$ for some even integer $n$. The ring $\mathbb{Z}_N$ will be identified with the interval $\{-K/2, -K/2 + 1, \ldots, N - K/2 - 1\}$ in the obvious way. For $0 \le i \le n + 1$, we call the set $[iK - K/2, iK + K/2)$ the *block $i$*.

During the procedure, we will label each particle as a *free* or a *carpet* particle. Free particles will be further subdivided into *frozen* or *thawed*. At the beginning, we declare all particles at sites $iK$, $i = 0, 1, \ldots, n + 1$, to be free and thawed, and all the other ones to be carpet particles.

We choose one free thawed particle (if there is one) to be the *hot* particle. The criterion for picking the hot particle will be described later.

There are some special sites, called *holes*, which will be "moved" during the procedure. At the beginning, we define the holes to be at sites $iK$ for $i = 0, 1, \ldots, n + 1$.

A *defect* or *defective site* is a site that is not a hole and has no particle. When an active particle reaches a defect, we say that the defect is *fixed*. In this case, the site stops being defective.

With the assumptions and definitions above, the initial configuration $\eta$ will satisfy the following properties, which will be preserved by the toppling procedure that will be described below.

(P1) Each block $i$ has exactly one hole, which is located at some site in

$[iK, iK + a]$.

(P2) Every site contains a carpet particle, except for the holes and defective sites.

(P3) If there are defective sites in $[iK - K/2, iK + K/2)$, the hole is at $iK$.

(P4) Carpet particles between the hole and $iK + a$ are active.

(P5) Free particles are always active.

(P6) All free particles are at sites $iK$ or $iK + a$, except, possibly, the hot particle.

(P7) There is at most one frozen free particle per block.

(P8) There is a frozen free particle in block $i$ if and only if the hole is at position $iK + a$. In this case, the frozen free particle in this block is at position $iK + a$.

(P9) The hot particle is free and thawed.

These properties hold for the initial configuration either by the definitions given on the beginning of this section or vacuously.

*Remark* 1.2. Properties (P1), (P4), (P5), (P6), (P7), (P8) and (P9) have identical counterparts in [HRR23]. If we assume the configuration has no defects, (P3) is vacuously true and (P2) becomes the same as in [HRR23].

## 1.2.2 Attempted emission

We will describe in this section a procedure for an *attempted emission*. An attempted emission can end either in a *successful emission* or in a *failure*.

Suppose we choose the hot particle at block $i$. A *successful emission* (or just *emission*) occurs when one of the following conditions happen:

- when the hot particle reaches $(i + 1)K$ (resp. $(i - 1)K + a$) and there are no defects in block $i + 1$ (resp. $i - 1$);

- when the hot particle reaches a vacant site in block $i + 1$ (resp. $i - 1$) and there are defects in block $i + 1$ (resp. $i - 1$). We remember that, by (P2), vacant sites are either holes or defects.

In either case, block $i$ is called the *emitting* block and block $i+1$ (resp $i-1$) the *receiving* block. When an emission occurs, we say that the hot particle is *emitted*.

Let us describe the procedure for an attempted emission. First, we will choose the hot particle using the following criterion. At the beginning of an attempted emission, pick the smaller $i$, $1 \leq i \leq n$, such that there are no defects and there is at least one thawed free particle in block $i$. If there is no such $i$, we declare the procedure finished. After $i$ is chosen, we pick one of the free thawed particles of block $i$ to be the hot particle. Note that these free thawed particles must be in $iK$ or $iK + a$ by (P6) If there is a free thawed particle at $iK$, we choose the hot particle at $iK$; else, we choose the hot particle at $iK + a$.

By (P5) and (P9) we can always topple the site which contains the hot particle. We will only topple the site where the hot particle is.

We have the following cases.

- **There is no frozen particle at site** $iK + a$

  Topple the site where the hot particle is until it is is emitted or reaches the hole in $[iK, iK+a]$. If it is emitted, we declare the attempted emission to be finished.

  If the hot particle reaches the hole and sleeps, we turn it into a carpet particle, move the hole to the site immediately to the right, turn the carpet particle in the new location of the hole into a free particle, and choose this particle to be hot ((P4) ensures that the new hot particle is active). If the hole reaches $iK + a$, we stop the attempted emission and declare the free particle at $iK + a$ as frozen. If this happens, we say that the attempted emission ends in a *failure* and, if possible a new hot particle is chosen in the way described above.

  If the particle leaves the hole, there are two sub-cases:

  - the particle is emitted. In this case, we finish the attempted emission.

  - the particle makes an excursion and returns to the hole. In this case, we move the hole to the leftmost site visited by the hot particle in this excursion, turn the hot particle into a carpet one, and turn the particle at the new position of the hole into a free one. (Note that, if the hot particle made an excursion to the right of the hole, the leftmost site visited is the hole itself.)

- **There is a frozen free particle at $iK + a$:**

  In this case, by the combination of (P2), (P3) and (P8), there are no defects on $[iK - K/2, iK + K/2)$, and all sites in $[iK, iK + a)$ have a carpet particle.

  We topple the site where the hot particle is until it is emitted.

  If the particle visits every site in $[iK, iK + a]$ before being emitted, we move the hole to $iK$, turn the frozen particle at $iK + a$ into a carpet particle and turn the carpet particle at $iK$ into a free particle. This finishes the attempted emission.

Note that, by (P8) those two cases are the only possibilities.

See Figure 1.2 to an example of a piece of the ring after some attempted emissions have been made.

## 1.3 Properties preserved

In this section, we prove that the properties (P1)-(P9) of a configuration are preserved by the procedure described in Section 1.2.

**Proposition 1.3.** *At the end of each attempted emission, each one of properties (P1)-(P9) is preserved.*

*Proof.* This follows from the following reasons:

(P1) We never move the hole outside $[iK, iK + a]$, and every time we move the hole to a new location, the old location ceases to be a hole. Thus, there is only one hole at each block.

(P2) This happens because sites that are either defects or holes only cease to be defects and holes if a particle reaches them and is turned into a carpet particle; and with exception of the hole, we only topple sites with at least two particles.

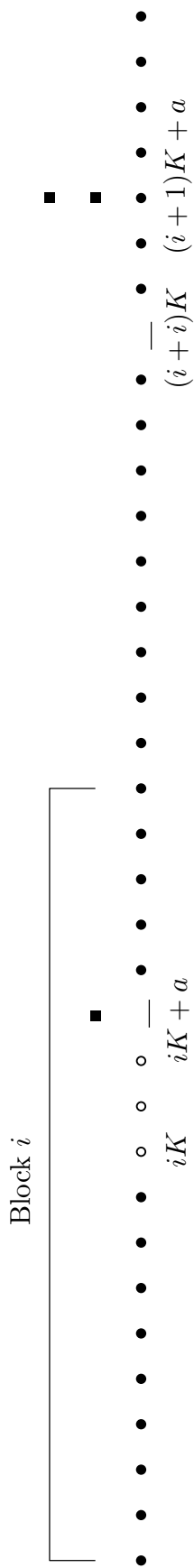

**Figure 1.2:** Example of possible state of a piece of the ring when an attempted emission is about to begin. Filled circles are active carpet particles, empty circles are sleeping carpet particles. The squares represent free particles, one of which is frozen. The holes are indicated by underlining the site.

(P3) We only move the hole after all defects in $[iK - K/2, iK + K/2)$ are fixed, because the block which the hot particle is chosen from has no defects.

(P4) If a particle fixes a defect, we relabel it as a carpet one and it remains active. If a particle sleeps, we move the hole to the right. We only move the hole to the left if all the particles between its new location and $iK + a$ are active.

(P5) Every time a free particle sleeps, we relabel it as a carpet particle.

(P6) We only move the hot particle, and every time the hot particle fixes a defect or sleeps, we relabel the hot particle as a carpet one. We only stop to move the hot particle when it is emitted or frozen. Thus, every time we choose a new hot particle, the old one is either relabeled as a carpet particle (if it fixed a defect), is frozen (at $iK + a$) or reached block $i - 1$ or $i + 1$.

(P7) If there is a frozen free particle at the block where the hot particle is chosen, the hot particle is always emitted or fixes a defect. In either case, we do not freeze any other particle.

(P8) This is because a particle is frozen if and only if the hole reaches $iK + a$, and is turned from frozen to thawed when the hole leaves $iK + a$.

(P9) We choose a free and thawed particle to be hot. When a hot particle is relabeled as a carpet one or declared frozen, a new hot particle is chosen.

$\square$

## 1.4   Alternating modes

We will now prove Theorem 1.1. To do so, we apply cyclically the procedure of Section 1.2 in two *modes*, A and B. At each mode, blocks 0 and $n+1$ are called *sinks* and blocks $n/2$ and $n/2 + 1$ are called *sources*. We apply the procedure of Section 1.2 until there are no possible choices for the hot particle outside the sink region. When this happens, we declare the mode to be ended. At the end of each mode we change the enumeration of the blocks in a way such that the sinks of mode A are the sources of mode B and vice-versa (see Figure 1.3). The idea behind of the proof of Theorem 1.1 is to show that, in each mode, if we start with many particles in the source, we can apply the procedure of Section 1.2 until there are many particles in the sink with high probability. Then, we will reverse the roles of the source of the sink and we can apply the procedure again. Thus, we conclude that, with high probability, we will repeat the procedure many times until the system is stabilized.

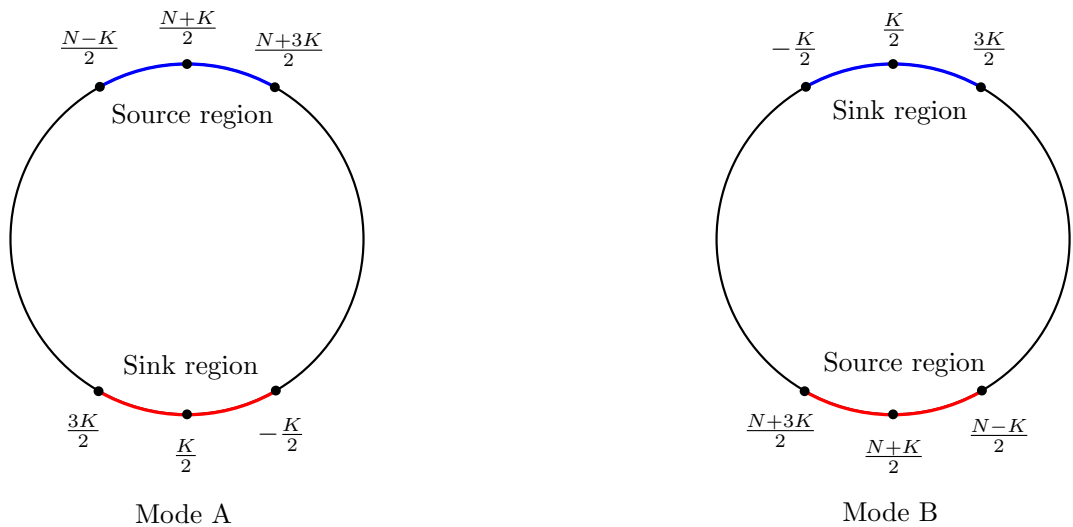

**Figure 1.3:** The two modes for the application of the dynamical procedure. At the end of each mode, the blocks are relabeled in a way that the source mode of each mode corresponds to the sink region of the next mode and vice-versa.

We remember that $n + 2 = N/K$ is the number of blocks.

Let $D$ be the number of defects at the end of a mode. We will define the following Condition and use it in the proof of Theorem 1.1.

**Condition 1.4.** *There are at least $\frac{7}{8}n + D$ free particles of which at most $\frac{5}{8}n$ are frozen.*

**Proposition 1.5.** *There is $c > 0$ depending on $\lambda$ such that, if at the beginning of a mode Condition 1.4 is satisfied, then at the beginning of the next mode Condition 1.4 will also be satisfied with probability at least $1 - e^{-cn}$, for all $n$ sufficiently large.*

We will prove Proposition 1.5 at the end of this section. Now we are ready to prove Theorem 1.1.

*Proof of Theorem 1.1.* We take $1 - \frac{1}{4K} < \zeta < 1$. Therefore, at the initial configuration there are at most $n/4$ defects, and at least $3n/4$ blocks with one free particle. Thus, there will be at least one block without defects and with one free thawed particle, and we can start the toppling procedure.

We will call a mode *successful* if Condition 1.4 is met at the end of the mode. We first note that, if a mode is successful, at least one jump is made, because there will be at least one block in the beginning of the mode with no defects and a free particle.

Using Proposition 1.5,

$$\mathbb{P}(\mathcal{J} \geq e^{\delta N}) \geq \mathbb{P}(\text{The first } e^{\delta N} \text{ modes are successful})$$
$$\geq (1 - e^{-cN/K})^{e^{\delta N}} \geq 1 - e^{(\delta - c/K)N}.$$

The theorem follows if we take $\delta < c/K$ and $\delta' = c/K - \delta$. $\qquad\square$

The following proposition states an important conservation law.

**Proposition 1.6.** *The number of free particles minus the number of defects is conserved during the procedure.*

*Proof.* There are two situations where a free particle is turned into a carpet one or vice versa. First, when a free particle fixes a defect. In this case, both the number of free particles and of defects decrease by one, so the difference is conserved. Second, when a free particle sleeps and it is turned into a carpet particle, the carpet particle immediately to its right is turned into a free one, and both numbers are conserved. This proves the proposition. $\qquad\square$

Let $E_n$ be the set of blocks that realize at least one successful emission during a mode, $F(E_n)$ be the number of frozen particles in $E_n$ and $F(\mathbb{Z}_N)$ the total number of frozen particles in the ring $\mathbb{Z}_N$.

To prove Proposition 1.5 we will make use of the following statements.

**Proposition 1.7.** *There exists a constant $c$ depending on $\lambda$ such that, if Condition 1.4 is satisfied then $\mathbb{P}(F(E_n) \geq \frac{n}{8}) \leq e^{-cn}$ for all $n$ sufficiently large.*

**Proposition 1.8.** *There exists a constant $c > 0$, depending on $\lambda$ such that, if at the beginning of the mode all the holes are at positions $iK$, then $\mathbb{P}(F(\mathbb{Z}_N) \geq \frac{n}{8}) \leq e^{-cn}$ for all $n$ sufficiently large.*

*Proof of Proposition 1.5.* As the number of free particles minus the number of defects is conserved by Proposition 1.6, at the end of the mode there are still $\frac{7}{8}n + D$ free particles.

In the first mode all the holes are at positions $iK$, $0 \leq i \leq n + 1$, hence by Proposition 1.8, the probability of having at most $n/8$ frozen particles in total at the ring is at least $1 - e^{-cn}$. In this event, Condition 1.4 is met.

Now, let us show that, on the event $\{F(E_n) \leq n/8\}$, Condition 1.4 is met at the end of the mode.

We first argue that in this event, there will be free particles inside the sink region. From the second mode onward, all free thawed particles begin the mode in the source region. Therefore, $E_n$ is a connected set surrounding the source region, and frozen particles outside $E_n$ must be in the two blocks adjacent to $E_n$.

If $\{F(E_n) \leq n/8\}$ happens, at the end of the mode, there are at most 2 new frozen particles outside $E_n$ and $5n/8$ particles that were frozen at the beginning of the mode. Since $5n/8 + n/8 + 2 < 7n/8$ (if $n > 16$), and the three options listed earlier are the only places where free frozen particles may be; there must have been free particles that are not frozen by the end of the mode, and this can only be achieved if free particles are emitted to the sink region.

Therefore, we conclude that, in fact, $E_n$ must contain the source region and be adjacent to the sink region, and since $E_n$ is a connected collection of blocks, the number of blocks in the complement of the union of $E_n$ and the sink region is at most $n/2$. In particular, by (P7), there are at most $n/2$ frozen particles outside $E_n$, hence at most $5n/8$ frozen particles in total.

By Proposition 1.7, $\mathbb{P}(F(E_n) < n/8) \geq 1 - e^{-cn}$, and Proposition 1.5 is proven. $\qquad \square$

# 1.5 Upper bound for the number of frozen particles

In this section, the analysis is restricted to the dynamics of a single mode. In Subsection 1.5.1 we will define some concepts that will be used in the proof of Proposition 1.7, In Subsection 1.5.2 we state a single-block estimate and use it to prove Proposition 1.7.

## 1.5.1 $\sigma$-algebras and flow of particles

To each site $y \in [iK - K/2, iK + K/2)$, $i = 0, i, \ldots, n+1$ we associate *three* independent copies of the stack of instructions, $\xi^y$, $\xi^{y,L}$ and $\xi^{y,R}$. Each time one of those sites is toppled we use the first unused instruction in $\xi^{y,L}$ if the hot particle is chosen at $\{(i-1)K, (i-1)K+a\}$, the first unused instruction in $\xi^y$ if the hot particle is chosen at $\{iK, iK+a\}$ or the first unused instruction in $\xi^{y,R}$ if the hot particle is chosen at $\{(i+1)K, (i+1)K+a\}$. Thus, each site $y$ in block $i$ has a stack associated to the block $i-1$ (except for $i = 0$), one stack associated to block $i$ and one stack associated to block $i-1$ (except for $i = n+1$).

We define the following $\sigma$-algebras

$$\begin{aligned}
\mathcal{F}_i = \sigma(\{\xi^y, \xi^{y,R}, \xi^{y,L} : y \in (-K/2, iK - K/2)\} \\
\cup \{\xi^{y,L}, \xi^y : y \in [iK - K/2, iK + K/2)\} \\
\cup \{\xi^{y,L} : y \in [iK + K/2, (i+1)K + K/2)\}).
\end{aligned}$$

In words, $\mathcal{F}_i$ comprises information about stacks associated with blocks $0, 1, \ldots, i$.

Let $\eta$ be the configuration of the ring at the beginning of a mode, and, for $0 \leq k \leq n$, define $\eta_k$ to be the restriction of $\eta$ to blocks $0, 1, \ldots, k$. Let $d'$ be the number of defects in $\eta_k$. We define, for $(m, d) \in \{0\} \times \{0, 1, \ldots, d'\} \cup \mathbb{N} \times \{d'\}$, $\eta'_k(m, d)$ in the following way: first, we use the toppling procedure of Section 1.2 restricted to blocks $0, 1, \ldots k$. After this is done, we define $\eta'_k(m, d)$ as the resulting configuration after block $k$ receives $m + d$ particles from block $k+1$ ($d$ particles for defective sites, $m-1$ particles arrive necessarily at $kK+a$ and 1 particle may arrive either at $kK$ if $kK$ is vacant and there are defects in

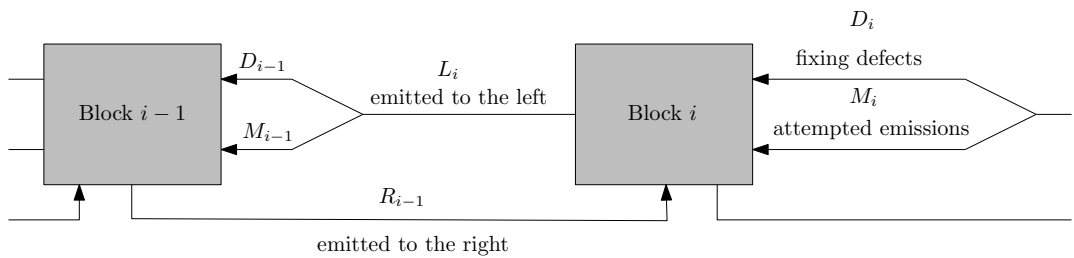

**Figure 1.4:** Pictorial representation of the mass-balance equation

$(kK - K/2, kK)$ or at $kK + a$ if not.)

Let $1 \leq i \leq k \leq n$, and $m$ and $d$ integers for which $\eta'_k(m, d)$ is well-defined. We run the procedure described on Section 1.2 on $\eta'_k(m, d)$ and define the following random functions:

- $L_i^k(m, d)$ is the number of particles that are emitted from block $i$ to block $i - 1$;

- $R_i^k(m, d)$ is the number of particles emitted from block $i$ to block $i + 1$;

- $D_i^k(m, d)$ is the number of particles that reached vacant sites in block $i$ from block $i + 1$ and will not be available for attempted emissions;

- $S_i^k(m, d)$ is the number of particles frozen in block $i$ after the procedure is finished.

By (P7), $S_i^k(m, d)$ can only be 0 or 1.

We also define $M_i^k$ as the number of particles emitted from block $i$ to block $i - 1$ that will be available for attempted emissions in block $i - 1$.

For the sake of simplicity, we may also write $L_i$, $R_i$, $D_i$, $M_i$ and $S_i$ instead of $L_i^i$, $R_i^i$, $D_i^i$, $M_i^i$ and $S_i^i$.

A particle that is emitted $i$ to block $i - 1$ can only fix a defect in block $i - 1$, reach the vacant hole or wait at $iK + a$ to be used in an attempted emission. Therefore, we can state the following mass-balance equation

$$L_i^k(M_i^k, D_i^k) = M_{i-1}^k + D_{i-1}^k. \tag{1.9}$$

The following proposition gives a relation between those random functions.

**Proposition 1.10.** *For every* $1 \leq i \leq n$,

$$S_i(M_i, D_i) = S_i^n(0, 0),$$
$$L_i(M_i, D_i) = L_i^n(0, 0),$$
$$R_i(M_i, D_i) = R_i^n(0, 0),$$
$$D_i^i(M_i, D_i) = D_i^n(0, 0).$$

*Proof.* For $i = n$, the result is immediate. Fix some $i \in \{1, \ldots, n-1\}$. Whenever possible, the procedure described on Section 1.2 chooses the hot particle in blocks $1, 2, \ldots, i-1$, and when this choice is possible, it depends only on the configuration in blocks $1, 2, \ldots, i-1$. When the choice in those blocks is not possible, the procedure either stops or waits for a particle to come to block $i$ from block $i+1$. This particle will either fix a defect in block $i$, reach the vacant hole at $iK$ or wait at $iK + a$ until it is chosen to be the hot particle. By the way the procedure was designed, these particles will have the least priority to be chosen as the hot particle, and will not affect the dynamics in the first $i$ blocks until they are chosen to be hot.

The only other way that attempted emissions in block $i+1$ can affect the outcome on blocks 1 to $i$ is by using instructions on stacks of instructions in sites between $iK + a$ and $(i+1)K$, but this is fixed by using multiple i.i.d. stacks of instructions in those sites. $\square$

### 1.5.2 Bounding the exponential expectation

To prove Proposition 1.7, we will use the following one-block estimate, to be proved in Section 1.6.

**Proposition 1.11.** *(i) For any initial configuration $\eta_0$ satisfying Properties (P1) - (P9),*

$$\sup_{\ell \geq 0, d \geq 0} \mathbb{E}\left[ \sum_{m=0}^{\infty} e^{48 S_i(m,d)} \mathbb{1}_{\{L_i(m,d) = \ell\}} \mathbb{1}_{\{L_i(m,d) + R_i(m,d) \geq 1\}} \,\Big|\, \mathcal{F}_{i-1} \right] \leq e^4.$$

*almost surely, where the sum runs over all $m$ for which $\eta_i'(m, d)$ is defined.*

*(ii) For any initial configuration $\eta_0$ satisfying Properties (P1) - (P9) and*

*with holes only at positions $iK$, $1 \le i \le n$,*

$$\sup_{\ell \ge 0, d \ge 0} \mathbb{E}\left[ \sum_{m=0}^{\infty} e^{48 S_i(m,d)} \mathbb{1}_{\{L_i(m,d) = \ell\}} \;\Big|\; \mathcal{F}_{i-1} \right] \le e^4,$$

*almost surely, where the sum runs over all $m$ for which $\eta_i'(m,d)$ is defined.*

*Proof of Proposition 1.7.* Denote, for simplicity, $G_{m,d}^i$ as the event $\{(R_i + L_i)(m,d) \ge 1\}$.

We will first show that

$$\mathbb{E}[e^{48 F(E_n)}] \le (n+1) 4^n e^{4n}.$$

By definition, $F(E_n)$ is the number of particles that are frozen on the collection of blocks that perform at least one emission. Also, at the beginning of the procedure there are at most $n$ free particles (one at each block), and the number of free particles never increases.

Since $\zeta > 1 - 1/4K$, we can bound the number of defects by $n$. Therefore, using the definition of $S_i$, Proposition 1.10 and the mass-balance equation (1.9):

$$\mathbb{E}[e^{48 F(E_n)}] =$$
$$\mathbb{E}\left[ \sum_{\substack{m_0, \ldots, m_n \\ d_0, \ldots, d_n}} \prod_{i=1}^{n} e^{48 S_i(m_i, d_i)} \mathbb{1}_{\{M_i^n = m_i\}} \mathbb{1}_{\{D_i^n = d_i\}} \mathbb{1}_{\{L_i(m_i, d_i) = m_{i-1} + d_{i-1}\}} \mathbb{1}_{G_{m_i, d_i}^i} \right]$$
$$\le \sum_{\substack{d_0, \ldots d_n \\ d_0 + \cdots + d_n \le n}} \mathbb{E}\left[ \sum_{m_0 = 0}^{n} \sum_{m_1, \ldots, m_n} \prod_{i=1}^{n} e^{48 S_i(m_i)} \mathbb{1}_{\{L_i(m_i) = m_{i-1} + d_{i-1}\}} \mathbb{1}_{G_{m_i, d_i}^i} \right].$$

The number of non-negative integer solutions of the inequality $d_0 + d_1 + \cdots + d_n \le n$ is bounded above by $4^n$. Thus,

$$\mathbb{E}[e^{48 F(E_n)}] \le 4^n \cdot \sup_{\{d_0, \ldots, d_n\}} \mathbb{E}\left[ \sum_{m_0 = 0}^{n} \sum_{m_1, \ldots, m_n} \prod_{i=1}^{n} e^{48 S_i(m_i)} \mathbb{1}_{\{L_i(m_i) = m_{i-1} + d_{i-1}\}} \mathbb{1}_{G_{m_i, d_i}^i} \right].$$

We will show by induction that, for every non-negative integer $r \le n$, and every possible choice of $d_0, \ldots d_n$,

$$\mathbb{E}\left[\sum_{m_0=0}^{n} \sum_{m_1,\ldots,m_r} \prod_{i=1}^{r} e^{48S_i(m_i,d_i)} \, \mathbb{1}_{\{L_i(m_i)=m_{i-1}+d_{i-1}\}} \mathbb{1}_{G^i_{m_i,d_i}}\right] \leq (n+1)e^{4r}.$$

Indeed, for $r = 0$ the expectation above is bounded by $n + 1$, and the base case of induction is finished. Suppose the inequality is valid for every integer between $0$ and $r - 1$. By conditioning on $\mathcal{F}_{r-1}$ and using item (i) of Proposition 1.11,

$$\mathbb{E}\left[\sum_{m_0=0}^{n} \cdots \sum_{m_r} \prod_{i=1}^{r} e^{48S_i(m_i,d_i)} \mathbb{1}_{\{L_i(m_i)=m_{i-1}+d_{i-1}\}} \mathbb{1}_{G^i_{m_i,d_i}}\right]$$

$$\leq \mathbb{E}\left[\sum_{m_0} \sum_{m_1} \cdots \sum_{m_{r-1}} \prod_{i=1}^{r-1} e^{48S_i(m_i,d_i)} \mathbb{1}_{\{L_i(m_i,d_i)=m_{i-1}+d_{i-1}\}} \mathbb{1}_{G^i_{m_i,d_i}}\right]$$

$$\cdot \sup_{\ell \geq 0, d \geq 0} \mathbb{E}\left[\sum_{m_r} e^{48S_r(m_r,d)} \mathbb{1}_{\{L_r(m_r,d)=\ell\}} \mathbb{1}_{G^r_{m_r,d_r}} \mid \mathcal{F}_{r-1}\right]$$

$$\leq e^4 \mathbb{E}\left[\sum_{m_0}^{n} \cdots \sum_{m_{r-1}} \prod_{i=1}^{r-1} e^{48S_i(m_i,d_i)} \mathbb{1}_{\{L_i(m_i,d_i)=m_{i-1}\}} \mathbb{1}_{G^i_{m_i,d_i}}\right]$$

$$\leq (n+1)e^4 e^{4(r-1)} = (n+1)e^{4r},$$

Therefore,

$$\mathbb{E}[e^{48F(E_n)}] \leq (n+1)4^n e^{4n}.$$

To conclude the bound, we note that, by Markov's inequality,

$$\mathbb{P}(F(E_n) \geq n/8) \leq e^{-6n} \mathbb{E}e^{48\mathrm{F}(E_n)}$$
$$\leq e^{-6n} \cdot (n+1)4^n e^{4n}$$
$$\leq e^{-cn}$$

for a constant $c > 0$ chosen properly and $n$ sufficiently large. $\qquad \square$

*Proof of Proposition 1.8.* The proof of this proposition is completely analogous to the proof of Proposition 1.7. First, we note that

$$\mathbb{E}[e^{48F(\mathbb{Z}_N)}] = \mathbb{E}\left[\sum_{\substack{m_0,\dots,m_n \\ d_0,\dots,d_n}} \prod_{i=1}^{n} e^{48S_i(m_i,d_i)} \mathbb{1}_{\{M_i^n=m_i\}} \mathbb{1}_{\{D_i^n=d_i\}} \mathbb{1}_{\{L_i(m_i,d_i)=m_{i-1}+d_{i-1}\}}\right]$$

$$\leq \sum_{\substack{d_1,\dots d_n \\ d_1+\cdots+d_n \leq n}} \mathbb{E}\left[\sum_{m_0=0}^{n} \sum_{m_1,\dots,m_n} \prod_{i=1}^{n} e^{48S_i(m_i)} \mathbb{1}_{\{L_i(m_i)=m_{i-1}+d_{i-1}\}}\right]$$

$$\leq 4^n \cdot \sup_{\{d_1,\dots,d_n\}} \cdot \mathbb{E}\left[\sum_{m_0=0}^{n} \sum_{m_1,\dots,m_n} \prod_{i=1}^{n} e^{48S_i(m_i)} \mathbb{1}_{\{L_i(m_i)=m_{i-1}+d_{i-1}\}}\right].$$

Since, in the first mode, all holes are at positions $iK$, item $(ii)$ of Proposition 1.11 can be used to prove by induction that, for every $r \leq n$,

$$\mathbb{E}\left[\sum_{m_0=0}^{n} \sum_{m_1,\dots,m_r} \prod_{i=1}^{r} e^{48S_i(m_i,d_i)} \mathbb{1}_{\{L_i(m_i)=m_{i-1}+d_{i-1}\}} \mathbb{1}_{G_{m_i,d_i}^i}\right] \leq (n+1)e^{4r},$$

and Markov's inequality is used as in the proof of Proposition 1.7. $\qquad \square$

## 1.6  One-block estimate

This section is devoted to proving Proposition 1.11.

We will begin by proving item $(i)$. Let $i$ be fixed. We will consider two cases. First, we will analyze the case where there are defects in $[iK-K/2, iK+K/2)$ after we stabilize $\eta_i'(m,d)$ using the procedure described in Section 1.2. The only value of $m$ compatible to the definition of $\eta_k'(m,d)$ in this case is $m=0$, and, by Properties (P3) and (P8), $S_i(0,d)=0$. It follows that

$$\sup_{\ell \geq 0, d \geq 0} \mathbb{E}\left[\sum_{m=0}^{\infty} e^{48S_i(m,d)} \mathbb{1}_{\{L_i(m,d)=\ell\}} \mathbb{1}_{G_{m,d}^i} \mid \mathcal{F}_{i-1}\right]$$
$$= \sup_{\ell \geq 0, d \geq 0} \mathbb{E}[e^{48S_i(0,d)} \mathbb{1}_{\{L_i(0,d)=\ell\}} \mathbb{1}_{G_{0,d}^i} \mid \mathcal{F}_{i-1}] \leq 1.$$

We now analyze the remaining case, i.e., when there are no defects in $[iK-K/2, (i+1)K/2)$ after we stabilize $\eta_i'(m,d)$ using the procedure described in Section 1.2. Thus, we need only to analyze the case where $d=0$.

We will use from now on a slight abuse of notation and denote by $L_i(j)$, $S_i(j)$ and $\mathbb{1}_{G_j^i}$ the values of the quantities $L_i$, $S_i$ and $\mathbb{1}_{G^i}$ after $j$ attempted emissions in block $i$.

For each particle added at $iK+a$, at least one attempted emission will be made. In other words, for each $m$ there is at least one corresponding $j$. Therefore,

$$\sum_{m=0}^{\infty} e^{48S_i(m,d)} \mathbb{1}_{\{L_i(m,d)=\ell\}} \mathbb{1}_{G_{m,d}^i} \leq \sum_{j=0}^{\infty} e^{48S_i(j)} \mathbb{1}_{\{L_i(j)=\ell\}} \mathbb{1}_{G_j^i}.$$

and it is sufficient to bound

$$\sup_{\ell \geq 0} \mathbb{E}\left[ \sum_{j=0}^{\infty} e^{48S_i(j)} \mathbb{1}_{\{L_i(j)=\ell\}} \mathbb{1}_{G_j^i} \mid \mathcal{F}_{i-1} \right] \tag{1.12}$$

For simplicity, as $i$ is fixed, we will write $S$, $R$, $L$ and $G$ instead of $S_i$, $R_i$, $L_i$ and $G^i$.

We will denote by $\mathcal{G}_j$ the $\sigma$-field generated by the information revealed up to the $j$-th attempted emission plus the information contained in $\mathcal{F}_i$.

**Definition 1.13.** The process $(H(j))_{j \in \mathbb{N}}$, taking values in $\{0, 1, \ldots, a\}$, is defined as $H(j) = v$ if the hole is at position $iK + v$ after the $j$-th attempted emission.

It is immediate that $H(j)$ is $\mathcal{G}_j$-measurable.

For brevity we will define $\tilde{\mathbb{P}}[\cdot] := \mathbb{P}[\cdot \mid \mathcal{F}_{i-1}]$ and $\tilde{\mathbb{E}}[\cdot] := \mathbb{E}[\cdot \mid \mathcal{F}_{i-1}]$.

To prove Proposition 1.11 we will make use of the following lemmas.

**Lemma 1.14.** *If $K$ is big enough, then, for every $j \geq 1$,*

$$\tilde{\mathbb{P}}[H(j) > a/2 \mid \mathcal{G}_{j-1}] \mathbb{1}_{\{H(j-1) \in [0,a/2] \cup \{a\}\}} \leq e^{-150}. \tag{1.15}$$

*and*

$$\tilde{\mathbb{P}}(a/2 < H(j) < a \mid \mathcal{G}_{j-1}) < e^{-150}. \tag{1.16}$$

*Remark* 1.17. The previous lemma is trivially true for $j = 0$ if the initial configuration has $H(0) = 0$.

**Lemma 1.18.** *For every $j \geq 0$, if $K$ is big enough then the following bound holds almost surely:*

$$\tilde{\mathbb{P}}(L(j+2) > L(j) \mid \mathcal{G}_j) \geq \frac{1}{5}.$$

For every integer $\ell \geq 0$, denote by $\tau_\ell$ the first $j$ for which $L(j) = \ell$, with $\tau_{-1} := -1$. We have that $\tau_{\ell+1} > \tau_\ell$, because every attempted emission will result in at most one particle emitted to the left. We can rewrite (1.12) to obtain

$$\mathbb{E}\bigg[\sum_{j=\tau_{i-1}+1}^{\infty} e^{48S(j)} \mathbb{1}_{\{L(j)=\ell\}} \mathbb{1}_{G_j} \mid \mathcal{F}_{i-1}\bigg]$$

$$\leq \sum_{j=\tau_{i-1}+1}^{\infty} \tilde{\mathbb{E}}[\mathbb{1}_{\{L(j)=\ell\}} \mathbb{1}_{\{S(j)=0\}} + e^{48} \mathbb{1}_{\{L(j)=\ell\}} \mathbb{1}_{G_j} \mathbb{1}_{\{S(j)=1\}}]$$

$$\leq \tilde{\mathbb{E}}[\tau_{\ell+1} - \tau_\ell] + \sum_{k=1}^{\infty} \tilde{\mathbb{E}}[e^{48} \mathbb{1}_{\{L(\tau_{\ell-1}+k)=\ell\}} \mathbb{1}_{G_{\tau_{\ell-1}+k}} \mathbb{1}_{\{S(\tau_{\ell-1}+k)=1\}}]. \quad (1.19)$$

Using Lemma 1.18, we get

$$\tilde{\mathbb{E}}[\tau_{\ell+1} - \tau_\ell] = \sum_{r=0}^{\infty} \tilde{\mathbb{P}}(\tau_{\ell+1} - \tau_\ell > r) \leq \sum_{r=0}^{\infty} \left(\frac{4}{5}\right)^{\lfloor r/2 \rfloor} < 10. \quad (1.20)$$

We split the sum in (1.19) at $k_0 := 250$. By Lemma 1.18 and our choice of $k_0$, we have that

$$\sum_{k=k_0}^{\infty} \tilde{\mathbb{E}}[e^{48} \mathbb{1}_{\{L(\tau_{\ell-1}+k)=\ell\}} \mathbb{1}_{G_{\tau_{\ell-1}+k}} \mathbb{1}_{\{S(\tau_{\ell-1}+k)=1\}}]$$

$$\leq \sum_{k=k_0}^{\infty} e^{48} \tilde{\mathbb{P}}(L(\tau_{\ell-1}+k)=\ell) \leq \sum_{k=k_0}^{+\infty} e^{48} \left(\frac{4}{5}\right)^{\lfloor k/2 \rfloor} < 1. \quad (1.21)$$

We will now prove that, for every $\ell \geq 0$ and $1 \leq k \leq k_0$,

$$\tilde{\mathbb{E}}[\mathbb{1}_{\{L(\tau_{\ell-1}+k)=\ell\}} \mathbb{1}_{G_{\tau_{\ell-1}+k}} \mathbb{1}_{\{S(\tau_{\ell-1}+k)=1\}}] < \frac{1}{k_0 e^{48}}. \quad (1.22)$$

Let us analyze first the case $\ell \geq 1$. In this case, obviously $\mathbb{1}_{G_{\tau_{\ell-1}}} = 1$. First, we consider $k = 1$. Note that if $S(\tau_{\ell-1} + 1) = 1$ and $L(\tau_{\ell-1} + 1) = \ell$, then the $(\tau_{\ell-1} + 1)$-th attempted emission is successful and $H(\tau_{\ell-1}) = a$. Using (1.15), after conditioning on $\mathcal{G}_{\ell-1}$,

$$\tilde{\mathbb{E}}[\mathbb{1}_{\{L(\tau_{\ell-1}+1)=\ell\}} \mathbb{1}_{\{S(\tau_{\ell-1}+1)=1\}}]$$

$$\leq \tilde{\mathbb{E}}[\mathbb{1}_{\{L(\tau_{\ell-1}+1)=\ell\}} \mathbb{1}_{\{H(\tau_{\ell-1}+1)=a\}} \mathbb{1}_{\{H(\tau_{\ell-1})=a\}}]$$

$$< e^{-150} < \frac{1}{k_0 e^{48}},$$

where the last inequality comes from our choice of $k_0$.

Still assuming $\ell \geq 1$, for $2 \leq k \leq k_0$, we have that

$$
\begin{aligned}
\tilde{\mathbb{P}}(S(\tau_{\ell-1} + k) = 0) &\geq \tilde{\mathbb{E}}[\mathbb{E}[\mathbb{1}_{\{H(\tau_{\ell-1}+k) \in [0, a/2]\}} \mid \mathcal{G}_{\tau_{\ell-1}+k-1}] \\
&\geq \tilde{\mathbb{E}}[(1 - e^{-150}) \mathbb{1}_{\{H(\tau_{\ell-1}+k-1) \in [0, a/2]\}}] \\
&= (1 - e^{-150}) \tilde{\mathbb{E}}[\mathbb{E}[\mathbb{1}_{\{H(\tau_{\ell-1}+k-1) \in [0, a/2]\}} \mid \mathcal{G}_{\tau_{\ell-1}+k-2}]] \\
&\geq (1 - e^{-150})^2 > 1 - \frac{1}{k_0 e^{48}},
\end{aligned}
$$

where in the second line we used (1.14) and in the last line we used (1.16).

This concludes the case $\ell \geq 1$. Now let us analyze the case $\ell = 0$. We remember that we defined $\tau_{-1} = -1$. When $k = 1$, $\mathbb{1}_{G_0} = 0$ trivially and (1.22) is valid. For $\ell = 0$ and $k = 2$, if the first attempted emission was not successful, then $\mathbb{1}_{G_1} = 0$. The only way the first attempt is successful, $S(1) = 1$ and $L(1) = 0$ occurs if $S(0) = 1$. By (1.15), the probability of $S(0) = 1$ and $S(1) = 1$ is less than $e^{-150}$. Either way, we conclude that (1.22) holds when $k = 1$.

For $\ell = 0$ and $2 < k \leq k_0$ we use (1.15) and (1.16) in the same way as before to conclude that

$$
\tilde{\mathbb{P}}(S(k) = 0) \geq (1 - e^{-150})(1 - e^{-150}) > 1 - \frac{1}{k_0 e^{48}}.
$$

This proves (1.22) and, therefore,

$$
\sum_{k=1}^{k_0} \tilde{\mathbb{E}}[e^{48} \mathbb{1}_{\{L(\tau_{\ell-1}+k) = \ell\}} \mathbb{1}_{G_{\tau_{\ell-1}+k}} \mathbb{1}_{\{S(\tau_{\ell-1}+k) = 1\}}] < 1. \tag{1.23}
$$

Combining (1.20), (1.21) and (1.23), we get

$$
\sup_{\ell \geq 0} \mathbb{E}\left[ \sum_{j=0}^{\infty} e^{48 S_i(j)} \mathbb{1}_{\{L_i(j) = \ell\}} \mathbb{1}_{G_j^i} \mid \mathcal{F}_{i-1} \right] < 10 + 1 + 1 < e^4 \tag{1.24}
$$

and the proof of case (i) Proposition 1.11 is complete.

The proof of $(ii)$ for $\ell \geq 1$ is almost identical to the proof of $(i)$. In the case $\ell = 0$, we cannot use Lemma 1.14 directly. But, since by hypothesis the initial configuration has holes only at $iK$, we can use Remark 1.17 instead of

Lemma 1.14 to show that, for $1 \le k \le k_0$,

$$\tilde{\mathbb{E}}\big[\mathbb{1}_{\{L(\tau_{\ell-1}+k)=\ell\}}\mathbb{1}_{\{S(\tau_{\ell-1}+k)=1\}}\big] < \frac{1}{k_0 e^{48}}$$

and proceed as in $(i)$.

## 1.7 Estimating the hole drift

This section is devoted to prove Lemmas 1.14 and 1.18.

Let $Z$ be a random variable measuring the maximum distance reached by a simple symmetric random walk before it returns to the origin, i.e., the law of $Z$ is given by $\mathbb{P}(Z = z) = \frac{1}{z(z+1)}$ for every positive integer $z$.

For each $v \in \mathbb{N}$, $v \le a$, define

$$Y_v = \begin{cases} +1, & \text{with probability } \frac{\lambda}{1+\lambda}; \\ 0, & \text{with probability } \frac{1/2}{1+\lambda}; \\ -\min(Z, v), & \text{with probability } \frac{1/2}{1+\lambda}. \end{cases}$$

We will use a modification of the variables $Y_v$ to obtain stochastic bounds for the hole drift when there is no emission. The probability of occurrence of an emission to the left (or to the right) is bounded above by

$$\delta := \frac{1}{K(1+\lambda)}.$$

Define

$$\tilde{Y}_v = \begin{cases} +1, & \text{with probability } \frac{\lambda}{1+\lambda}; \\ 0, & \text{with probability } \frac{1/2}{1+\lambda} + \delta; \\ -k, & \text{with probability } \frac{1/2}{1+\lambda}\frac{1}{k(k+1)} \text{ for } k = 1, 2, \dots, v-1; \\ -v & \text{with probability } \frac{1/2}{1+\lambda} - \delta - \sum_{k=1}^{v-1}\frac{1/2}{1+\lambda}\frac{1}{k(k+1)}. \end{cases}$$

Those random variables dominate stochastically the change of the position of the hole, when it is at $iK + v$, conditioned to the event that no emission is made and there are no defects inside the block.

**Lemma 1.25.** *If $a$ is large enough and $v \geq a/3$,*

$$\mathbb{E}[Y_v], \mathbb{E}[\tilde{Y}_v] \leq -40.$$

*Proof.* Note that

$$\mathbb{E}[Y_v] \leq \frac{\lambda}{\lambda+1} - \frac{1}{2(\lambda+1)} \sum_{k=1}^{a/3} \frac{1}{k(k+1)} \cdot k$$
$$\leq \frac{\lambda}{\lambda+1} - \frac{\log(a/3) - \log 2}{2(\lambda+1)},$$

and

$$\mathbb{E}[\tilde{Y}_v] = \mathbb{E}[Y_v] + v\delta.$$

The lemma follows if we take $a$ greater than $6 \exp(120(\lambda+1))$. $\qquad \square$

**Lemma 1.26.** *Let $(\tilde{Y}_{a/3}^{(i)})_{i \in \mathbb{N}}$ be a sequence of i.i.d. random variables with common distribution $\tilde{Y}_{a/3}$. Then, there is a constant $\alpha > 0$ such that*

$$\mathbb{P}\left(\sum_{i=1}^{a/6} \tilde{Y}_{a/3}^{(i)} \geq -a/6\right) \leq e^{-\alpha a}$$

*and*

$$\mathbb{P}\left(\sum_{i=1}^{a/2} \tilde{Y}_{a/3}^{(i)} \geq -2a/3\right) \leq e^{-\alpha a}$$

*Proof.* We will use the following Chernoff bound. Let $(Y^i)_{i=1}^{+\infty}$ be a sequence of bounded i.i.d. random variables with expectation $\nu$. Fix $b > 0$. For every $\gamma > 0$ such that $\gamma b$ is an integer with $\nu < -\gamma^{-1}$, there is a constant $\alpha' > 0$ such that

$$\mathbb{P}\left(\sum_{i=0}^{\gamma b} Y^i > -b\right) \leq e^{-\alpha' \gamma b}.$$

By Lemma 1.25, $\mathbb{E}[\tilde{Y}_{a/3}] < -40$. We can then take $b = -a/6$ and $\gamma = 1$ to obtain the first inequality of the lemma, and $b = -2a/3$, $\gamma = 3/4$ to obtain the second one. $\qquad \square$

We call a *step* of the toppling procedure when the hot particle starts in the hole and either

(i) falls asleep;

(ii) is emitted; or

(iii) leaves the hole and returns to it.

Let $T_j$ be the number of steps taken by the procedure between the $(j-1)$-st and the $j$-th attempted emissions.

**Lemma 1.27.** *For all $v \in \mathbb{Z}_N \cap [0, a]$,*

$$\mathbb{P}(T_j > a^3 \mid H(j-1) = v) < \frac{1}{a},$$

*and, if $v < a/2$,*

$$\mathbb{P}(T_j < a/2 \mid H(j-1) = v) < \frac{1}{2a}.$$

*Proof.* First we note that if the hot particle is at the hole, the probability that the next instruction will be a jump instruction is $\frac{1}{\lambda+1}$. Therefore, the probability of an emission is at least $\frac{1}{(\lambda+1)(3K/2+a)}$.

The first inequality comes then from the fact that

$$\left(1 - \frac{1}{(\lambda+1)(3K/2+a)}\right)^{a^3} < \frac{1}{a},$$

if $a$ is big enough.

If $v < a/2$, then the first $a/2$ steps of an attempted emission cannot end with a failure. If $T_j < a/2$, then the $j$th attempted emission ended in an successful emission. The probability of occurrence of an emission is at most $\frac{1}{(\lambda+1)K}$, then an union bound gives that the probability of an emission occurring in the first $a/2$ steps is at most

$$(a/2)\frac{1}{(\lambda+1)K} \le \frac{a}{2K} = \frac{1}{2a}$$

and the second inequality of the lemma follows. $\square$

*Proof of Lemma 1.14.* We begin by proving (1.15). As $\mathcal{G}_{j-1}$ contains no information after the $(j-1)$-st attempted emission, it is sufficient to prove that, for every $v \in ([0, a/2] \cup \{a\}) \cap \mathbb{Z}_N$,

$$\mathbb{P}(H(j) > a/2 \mid H(j-1) = v) < \frac{4}{a} < e^{-150}.$$

We split in three cases.

- $v = a$. In this case, by (P8) every site in the block has a carpet particle, and the hot particle is performing a simple random walk. If the hot particle visits every site in $[0, a]$ before being emitted, the hole is reset to position 0. The probability of this occurring is at least

$$\frac{K/2 - a}{K/2} = 1 - \frac{2}{a}$$

as desired.

- $v \leq a/3$. By Lemma 1.27, it is immediate that

$$\tilde{\mathbb{P}}(\{H(j) > a/2\} \cap \{T_j \geq a^3\} \mid H(j-1) = v) < \frac{1}{a}.$$

If $\{H(j) > a/2\} \cap \{T_j < a^3\}$ occurs, one of the first $a^3$ times the hole is at $a/3$ it moves to $a/2$ before returning to the left of $a/3$. This probability is bounded above by the probability that the sum of $a/6$ independent copies of $\tilde{Y}_{a/3}$ is at least 0. Using Lemma 1.26,

$$\tilde{\mathbb{P}}(\{H(j) > a/2\} \cap \{T_j > a^3\} \mid H(j-1) = v) < \frac{1}{a}.$$

The case follows by an union bound.

- $a/3 < v \leq a/2$. A failure cannot occur in the first $a/2 - 1$ steps. The probability that an emission occurs and $T_j < a/2$ is bounded above by $1/4a$ by Lemma 1.27. If an emission does not occur in the first $a/2 - 1$ steps and the hole is always at the right of $a/3$, the movement of the hole is stochastically dominated by $\tilde{Y}_{a/3}$. Thus the probability of $T_j < a/2$ and the hole be always at the right of $a/3$ is bounded above by the probability that the sum of $a/6$ copies of $\tilde{Y}_{a/3}$ is at least $-a/6$, and this probability is at most $e^{-\alpha a}$ by Lemma 1.26. Once the hole is at or at the left of $a/3$ we proceed as in the previous case. Hence, we get

$$\tilde{\mathbb{P}}(\{H(j) > a/2\} \cap \{T_j < a/2\} \mid H(j-1) = v) < \frac{1}{a},$$

and

$$\tilde{\mathbb{P}}(\{H(j) > a/2\} \cap \{a/2 \le T_j < a^3\} \mid H(j-1) = v) < \frac{1}{a}.$$

By Lemma 1.27,

$$\tilde{\mathbb{P}}(\{H(j) > a/2\} \cap \{T_j > a^3\} \mid H(j-1) = v) < \frac{1}{a},$$

and we conclude by an union bound.

This finishes the proof of (1.15).

It remains to prove (1.16). Again, it is sufficient to condition on $H(j-1)$. If $H(j-1) \in [0, a/2] \cup \{a\}$, we can use (1.15) to bound the probability of $H(j) \in (a/2, a)$. So let us suppose $H(j-1) \in (a/2, a)$. Since when a failure occurs we have $H(j) = a$, we only need to consider the case of an emission. The probability that an emission occurs in the first $a/2$ steps is bounded above by $\frac{1}{2a}$ by the same argument used in the proof of Lemma 1.27. The probability that a emission does not occur and the hole is never in $[0, a/3]$ in the first $a/2$ steps is at most $e^{-\alpha a}$ by Lemma 1.26. Once the hole reaches $[0, a/3]$ we can use Lemma 1.26 again as in the proof of (1.15) to conclude that the probability the hole reaches $(a/2, a)$ is at most $e^{-\alpha a}$. An union bound finishes the proof. □

*Proof of Lemma 1.18.* We notice that if an attempted emission ends in a failure, the next hot particle will be emitted. Thus, for each two attempted emissions, at least one will be successful. The situation where is a emission to the left is least likely is when the hot particle is chosen at $iK + a$ and there is a defect in $(i-1)K - K/2 + 1$. The probability of this emission occurs to the left is at least the probability that a random walk starting at $iK + a$ reaches $iK - 3K/2$ before it reaches $(i+1)K$, i.e., at least $\frac{((i+1)K/2) - (iK+a)}{((i+1)K/2) - (iK - 3K/2 + 1)} = \frac{1}{4} - \frac{1}{\sqrt{K}} \ge \frac{1}{5}$. □

# Part II

# Equilibrium and quasi-equilibrium in interacting particle systems

# Chapter 2

# Convergence to equilibrium for the critical contact process with modified border

The Harris contact process is a widely studied model for the spread of infections. Each site of $\mathbb{Z}$ (that we may think as an individual in a population) may be *infected* or *healthy*. Infected individuals become healthy at rate 1 and at rate $2\lambda, \lambda > 0$, it chooses a neighbor at random to try to infect. If this neighbor is healthy, it becomes infected. We will assume that the reader is familiar with the basic concepts and results about the contact process as shown in [Gri81] or [Lig05, Chapter VI].

The contact process exhibits a critical behavior: if the infection rate $\lambda$ is too small, infection tends to die out; if it is big, infection tends to persist forever. In other words, suppose that at time $t = 0$ there is only one infected individual at the origin. There is a *critical value* $\lambda_c$ such that, if $\lambda < \lambda_c$ then the infection will eventually disappear a.s. (i.e., all individuals will be healthy); and if $\lambda > \lambda_c$, there is a positive probability that there will be always infected individuals. At criticality ($\lambda = \lambda_c$), the infection will die out a.s. as shown in [BG90].

In [Dur84], Durrett proved of existence of an invariant measure for the critical and supercritical classical contact process seen from the right edge, that is, the contact process translated in such a way that its rightmost infected particle is at the origin. Galves and Presutti proved convergence of the supercritical process contact seen from the edge to the invariant measure [GP87] and this result was extended to the critical case in [CDS91]. Nonexistence of a invariant measure for the subcritical process was proved in [Sch87] in the discrete time case, and in [ASS90] for the continuous time case. However, when conditioned on non-extinction, subcritical contact processes converge to a distribution supported on infinite sets (see [AEGR15] and Chapter 3 of this thesis).

Many modifications of the basic contact process have been proposed. One of them, described in [DS00], is the *contact process with modified border*. In this process, infections to the right of the rightmost site and to the left of the leftmost site occur at rate $\lambda_e$ and infections at other sites occur at rate $\lambda_i$.

Given a measure $\mu$ on the power set of $\mathbb{Z}$, we denote by $\xi^\mu = (\xi_t^\mu)_{t\geq 0}$ to be a contact process with parameters $(\lambda_i, \lambda_e)$ and initial condition sampled from $\mu$ and denote by $\Psi\xi^\mu$ the process $\xi^\mu$ as seen from the right edge. For a subset $A$ of $\mathbb{Z}$ we write $\xi^A$ and $\Psi\xi^A$ for $\xi^{\delta_A}$ and $\Psi\xi^{\delta_A}$ respectively.

We denote by $\theta(\lambda_i, \lambda_e)$ the probability that this process survives forever when at time 0 there is only one infected particle.

In the case of the modified border contact process, the following result was first proved in [AR23].

**Theorem 2.1.** *Suppose $\lambda_e \leq \lambda_i$ and $\theta(\lambda_i, \lambda_e) > 0$ or $\lambda_e = \lambda_i = \lambda_c$. There is a measure $\mu$ such that for every initial condition $A \subseteq \mathbb{Z}$ with $|A| = +\infty$ and $\sup A < +\infty$), $\Psi\xi_t^A \to \mu$ weakly. Moreover, $\Psi\xi_t^\mu \sim \tilde{\mu}$ for every $t \geq 0$.*

When $\lambda_e \leq \lambda_i$, the process is *attractive*. This means that if we add more particles to the initial configuration, it helps the infection to survive (see [Lig05, Chapter III] for a formal definition). This is a direct consequence of the construction of the process given in Section 2.1. In contrast, when $\lambda_i < \lambda_e$, the process is non-attractive. The tools and results available for the study of non-attractive processes are far more limited than the ones for the attractive case.

In this work, we extend Theorem 2.1 to the non-attractive case where $\lambda_e = \lambda_c + \varepsilon$, $\lambda_i = \lambda_c$.

**Theorem 2.2.** *Suppose $\lambda_e = \lambda_c + \varepsilon$ and $\lambda_i = \lambda_c$. There is a measure $\tilde{\mu}$ such that for every initial condition $A \subseteq \mathbb{Z}$ with $|A| = +\infty$ and $\sup A < +\infty$), $\Psi\xi_t^A \to \tilde{\mu}$ weakly. Moreover, $\Psi\xi_t^{\tilde{\mu}} \sim \mu$ for every $t \geq 0$.*

Let $\mathcal{R}\xi_t^{\tilde{\mu}}$ be the position of the right edge of the configuration $\xi_t^{\tilde{\mu}}$, that is, $\mathcal{R}\xi_t^{\tilde{\mu}} = \sup \xi_t^{\tilde{\mu}}$. Theorem 2.2 implies that $\mathcal{R}\xi_t^{\tilde{\mu}}$ has stationary increments. By Birkhoff's Ergodic Theorem, $n^{-1}\mathcal{R}\xi_n^{\tilde{\mu}}$ converges to a.s. to a random variable $V$. We can extend this convergence along the naturals to convergence along al $t$ by noting that, if we ignore all deaths, the increments of the right edge will be distributed like a Poisson process $N$ of parameter $\lambda_e$. Thus,

$$\mathbb{P}\left(\max_{n \le t \le n+1}(\mathcal{R}\xi_t^{\tilde{\mu}} - \mathcal{R}\xi_n^{\tilde{\mu}}) \ge \varepsilon n\right) \le \mathbb{P}(N \ge \varepsilon n),$$

$$\mathbb{P}\left(\max_{n-1 \le t \le n}(\mathcal{R}\xi_n^{\tilde{\mu}} - \mathcal{R}\xi_t^{\tilde{\mu}}) \ge \varepsilon n\right) \le \mathbb{P}(N \ge \varepsilon n)$$

and then apply the Borel-Cantelli lemma. Thus, we have the following corollary.

**Corollary 2.3.** *There is a random variable $V$ with $\mathbb{E}[V] = \mathbb{E}[\mathcal{R}\xi_1^{\tilde{\mu}}]$ such that*

$$\frac{\mathcal{R}\xi_t^{\tilde{\mu}}}{t} \xrightarrow{a.s.} V.$$

*Remark* 2.4. It still an open problem to prove that the above random variable is degenerate, i.e., to prove existence of speed of the edge process.

The technique used in [AR23] to prove invariance of the limiting distribution does not translate immediately to the critical case. Theorem 2.2 is proved with arguments similar to [CDS91, Theorem 1], using the fact that, the infections starting at two different infinite initial conditions will tend to agree on finite dimensional cylinder sets.

The characteristics of the phase diagram of the contact process with modified border begun to be investigated in [DS00]. Andjel and Rolla [AR23] established almost entirely the characteristics of the phase diagram, but left open the question of survival of the process in the critical line of the attractive region. We solve this question with the following result.

**Theorem 2.5.** *Define*

$$\lambda_*^e(\lambda_i) = \inf\{\lambda_e : \theta(\lambda_i, \lambda_e) > 0\}.$$

*Then, if $\lambda_i > \lambda_c$, $\theta(\lambda_i, \lambda_*^e(\lambda_i)) = 0$.*

By [AR23, Theorems 4 and 5], we have that $\lambda_*^e(\lambda_i) < \lambda_i$ when $\lambda_i > \lambda_c$. As the process is attractive if $\lambda_i \ge \lambda_e$, the proof of Theorem 2.5 follows the steps

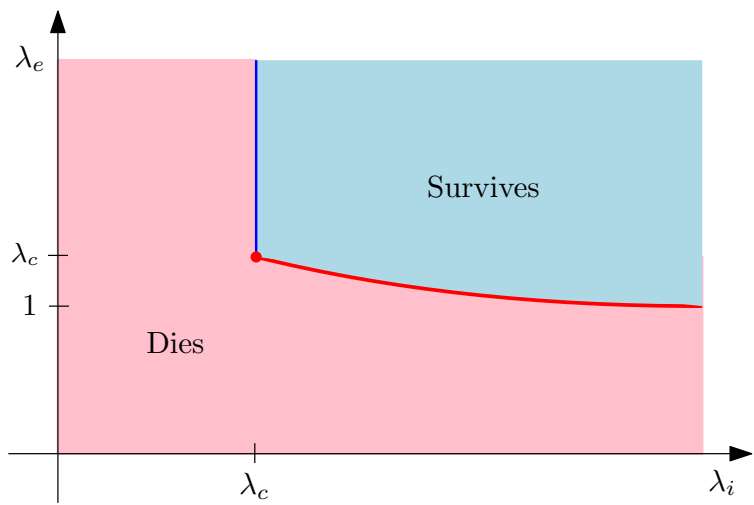

**Figure 2.1:** Phase space for the contact process with modified borde

of the proof of the analogous result for the usual critical contact process of [BG90].

With Theorem 2.5, the phase space of the contact process with modified border is fully characterized, as seen in Figure 2.1.

This chapter is structured as follows. Section 2.1 is dedicated to the construction of the contact process with modified border. Theorem 2.2 is proved in Section 2.2 and Theorem 2.5 is proved in Section 2.3.

## 2.1 Description of the model

A *configuration* is a subset of $\mathbb{Z}$. The set of all configurations is denoted $\Sigma$, i.e., $\Sigma$ is the power set of $\mathbb{Z}$, which we will identify with $\{0, 1\}^{\mathbb{Z}}$. We take the product topology in $\Sigma$, and thus $\Sigma$ is compact. For $A \subseteq \mathbb{Z}$, let $\mathcal{R}A = \sup A$ and $\mathcal{L}A = \inf A$. We define the right shift operator $T : \Sigma \to \Sigma$ such that $TA = A + 1 := \{x + 1 : x \in A\}$.

Let $\Sigma^{\ominus}$ be the set of semi-infinite configurations to the left, i.e., $\Sigma^{\ominus} = \{A \in \Sigma : |A| = +\infty, \mathcal{R}A < +\infty\}$, and $\Sigma^{\odot} = \{A \in \Sigma^{\ominus} : \mathcal{R}A = 0\}$. We also define the operator $\Psi$ in $\Sigma^{\ominus}$ as $\Psi A = T^{-\mathcal{R}A}A$.

We use an enlarged version of the usual construction of the contact process to couple processes with different parameters. We do this in the following way. For each $x$ in $\mathbb{Z}$ we sample a Poisson point process (PPP) $\omega_x$ with intensity 1.

We say there is a recovery mark at $(x,t)$ if $t \in \omega_x$. We sample also a collection of PPPs $\omega_{x,x\pm1} \subseteq (0,+\infty) \times (0,+\infty)$ with intensity 1 and say that the oriented edge from $(x,t)$ to $(x\pm1,t)$ is $\lambda$-open if $\omega_{x,x\pm1} \cap (\{t\} \times [0,\lambda)) \neq \emptyset$. We denote by **P** the probability measure on $\mathbb{Z} \times [0,+\infty)$ induced by those Poisson processes.

Take $A \subseteq \mathbb{Z}$ finite. We construct the contact process with modified border $\eta_{s,t}^A$ as a càdlàg process on $\Sigma$ with parameters $(\lambda_i, \lambda_e)$ as follows. Define $\eta_{s,s}^A = A$. Let $B$ be the state of the process at time $t^-$. Given $x \in B$, at time $t$, if there is a recovery mark at $(x,t)$ the process jumps from $B$ to $B \setminus \{x\}$; and if $x = \mathcal{L}B$ and the edge from $(x,t)$ to $(x-1,t)$ is $\lambda_e$-open or $x \neq \mathcal{L}B$ and the edge from $(x,t)$ to $(x-1,t)$ is $\lambda_i$-open, the process jumps from $B$ to $B \cup \{x-1\}$. Similarly, f $x = \mathcal{R}B$ and the edge from $(x,t)$ to $(x+1,t)$ is $\lambda_e$-open or $x \neq \mathcal{R}B$ and the edge from $(x,t)$ to $(x+1,t)$ is $\lambda_i$-open, the process jumps at time $t$ from $B$ to $B \cup \{x+1\}$. For $A$ infinite we define $\eta_{s,t}^A$ by taking a limit of $\eta_{s,t}^{A_n}$ with $A_n$ finite, $A_n \nearrow A$. If $s = 0$, we can write $\eta_t^A$ for $\eta_{0,t}^A$. If $A = \mathbb{Z}^-$, the set of non-positive integers, we write $\eta_t^-$ for $\eta_t^A$. We denote $\mathbb{P}_{\lambda_i,\lambda_e}$ the law of the contact process with parameters $(\lambda_i, \lambda_e)$. The subscripts may be dropped if the parameters are clear from context.

A function $f : \Sigma \to \mathbb{R}$ is said to be *increasing* if $f(A) \leq f(B)$ for every $A \subseteq B$. When $\xi$ and $\eta$ are random elements of $\Sigma$, we say that $\xi$ *stochastically dominates* $\eta$, and denote it by $\xi \succcurlyeq \eta$ if $\mathbb{E}[f(\xi)] \geq \mathbb{E}[f(\eta)]$ for every bounded increasing function $f$.

A *path* is a cadlag function $\gamma : [t_0, t_1] \to \mathbb{Z}$ whose jumps are of length 1. A path $\gamma$ defined in $[t_0, t_1]$ is said to be $(\lambda_i, \lambda_e)$-*active* for a configuration $A$ if $\gamma_s \in \eta_{t_0,s}^A$ for all $s \in [t_0, t_1]$ and for every point of discontinuity $u$ of $\gamma$ there is a $\lambda_i$ or $\lambda_e$- open edge from $\gamma_{u^-}$ to $\gamma_u$ .

## 2.2 Convergence to invariant measure

This section is devoted to the proof of Theorem 2.2.

*Proof of Theorem 2.2.* By compactness of $\Sigma$, there is a sequence $(t_k)_k$ converging to infinity, such that $\Psi \xi_{t_k}^A$ converges weakly to a measure $\tilde{\mu}$. It remains to prove that $\tilde{\mu}$ does not depend on the initial condition $A$ and that convergence holds for every $t$, and not just along that subsequence. Those are consequences of the following proposition.

**Proposition 2.6.** *For every $A$ and $B$ in $\Sigma^\odot$ and every finite dimensional cylinder set $C$,*

$$|\mathbb{P}(\Psi\xi_t^A \in C) - \mathbb{P}(\Psi\xi_t^B \in C)| \to 0,$$

*uniformly over $A$ and $B$.*

To prove that convergence holds not just along the subsequence $t_k$, but for every $t$, we note that for every $A \in \Sigma^\odot$, denoting by $\mu_{t-t_k}$ the law of the process at time $t - t_k$,

$$|\tilde{\mu}(C) - \mathbb{P}(\Psi\xi_t^A \in C)| \leq |\tilde{\mu}(C) - \mathbb{P}(\Psi\xi_{t_k}^A \in C)| + |\mathbb{P}(\Psi\xi_{t_k}^A \in C) - \mathbb{P}(\Psi\xi_t^A \in C)|$$
$$= |\tilde{\mu}(C) - \mathbb{P}(\Psi\xi_{t_k}^A \in C)|$$
$$+ \left| \int [\mathbb{P}(\Psi\xi_{t_k}^A \in C) - \mathbb{P}(\Psi\xi_{t_k}^B \in C)]\mathrm{d}\mu_{t-t_k}(B) \right|.$$

Using convergence of $\Psi\xi_t^A$ to $\tilde{\mu}$ along the subsequence $(t_k)_k$ and uniform convergence given by Proposition 2.6, we conclude that, in fact, $\Psi\xi_t^A$ converges weakly to $\tilde{\mu}$.

It remains to show that the limiting measure does not depend on the initial condition. Indeed, if $\Psi\xi_t^A \to \tilde{\mu}$ and $\Psi\xi_t^B \to \tilde{\nu}$, then by Proposition 2.6 $\tilde{\mu}$ and $\tilde{\nu}$ must agree on finite dimensional cylinders, implying $\tilde{\mu} = \tilde{\nu}$.

To prove invariance of $\tilde{\mu}$, we make use of the following proposition.

**Proposition 2.7.** *The measure $\tilde{\mu}$ satisfies $\tilde{\mu}(\Sigma^\odot) = 1$.*

For every fixed $t \geq 0$ and every continuous and bounded $f : \Sigma \to \mathbb{R}$ the functional $\mathcal{S}_t f$ given by $A \mapsto \mathbb{E}[f(\Psi\xi_t^A)]$ is continuous at every $A \in \Sigma^\odot$. Then, by weak convergence of $\Psi\xi_t^A$ to $\tilde{\mu}$. This means that $(\xi_t^{\tilde{\mu}})_t$ is a Feller process on $\Sigma^\odot$, and by Proposition 2.7 and [Lig05, Proposition I.1.8], it follows that $\xi_t^{\tilde{\mu}} \sim \tilde{\mu}$ for all $t > 0$. $\qquad\qquad\square$

The proof of Proposition 2.6 is a reformulation of the argument of [CDS91] with some modifications. The main difference arises from the fact that, in the classical contact process, we can couple the processes $\xi_t^0$ and $\xi_t^-$ such that $\xi_t^0 \cap [-L, 0] = \xi_t^- \cap [-L, 0]$ if $t$ is large enough, for all $L > 0$. This is no longer true in our case, but this issue can be easily solved if we consider a process with infection $\lambda_c + \varepsilon$ only at the right border and $\lambda_c$ everywhere else. This process also survives (see the first Remark after Theorem 5 of [AR23]).

*Proof of Proposition 2.6.* Take $(\zeta_t)_t$ to be a contact process that has parameters $\lambda_c + \varepsilon$ on the right edge and $\lambda_c$ everywhere else. That means that infections occurs from $\mathcal{R}\zeta_t$ to $\mathcal{R}\zeta_t + 1$ at rate $\lambda_c + \varepsilon$, from $\mathcal{R}\zeta_t$ to $\mathcal{R}\zeta_t - 1$ with rate $\lambda_i$ and occurs from every other infected site (including $\mathcal{L}\zeta_t$) with rate $\lambda_i$.

Define a sequence of stopping times $(\tau_n)_n \geq 0$ by $\tau_0 := 0$ and $\tau_k := \inf\{t > \tau_{k-1} : |\zeta^0_{\tau_{k-1},t}| = 0\}$. Take $A \in \Sigma^{\odot}$. We will define the process $\tilde{\zeta}^A_t$ for $t \in [\tau_{k-1}, \tau_k)$ as

$$\tilde{\zeta}^A_t = \zeta^{\Psi\eta}_{\tau_{k-1},t} + \mathcal{R}\eta,$$

with $\eta = \tilde{\zeta}^A_{\tau_{k-1}}$. By translation invariance of the construction of the contact process described in Section 2.1, $\Psi\tilde{\zeta}^A$ has the same distribution of of $\Psi\zeta^A$. Define also $N(t) = \sum_{k=1}^{\infty} \mathbb{1}_{\{\tau_k \leq t\}}$.

Let $L < +\infty$ fixed, $\kappa \in \Sigma^{\odot}$ and consider the cylinder $C = \{D \in \Sigma^{\odot} : D \cap [-L, 0] = \kappa \cap [-L, 0]\}$. We note that every cylinder may be represented in this form for some $L$ and $\kappa$.

For every $A, B \in \Sigma^{\odot}$,

$$|\mathbb{P}(\Psi\xi^A_t \in C) - \mathbb{P}(\Psi\xi^B_t \in C)| = |\mathbb{P}(\Psi\tilde{\zeta}^A_t \in C) - \mathbb{P}(\Psi\tilde{\zeta}^B_t \in C)|$$
$$\leq \mathbb{P}(\Psi\tilde{\zeta}^A_t \cap [-L, 0] \neq \Psi\tilde{\zeta}^B_t \cap [-L, 0]).$$

By the coupling above, the last probability is no larger than

$$\mathbb{P}(|\zeta^0_{\tau_{N(t)},t}| \leq L).$$

As $\mathbb{P}(|\zeta^0_s| \leq L \mid \tau_0 > s) \to 0$ as $s \to +\infty$ and by [AR23, Theorem 5] there is a positive probability that $\zeta^0_t \neq \emptyset$ for all $t$. Hence, $t - \tau_{N(t)} \to +\infty$ almost surely, so

$$|\mathbb{P}(\Psi\xi^A_t \in C) - \mathbb{P}(\Psi\xi^B_t \in C)| \to 0.$$

$\square$

The proof of Proposition 2.7 is based on [CDS91], that uses some results of [ASS90]. However, in the non-attractive regime, there is no deterministic velocity of the right edge, and the argument of [ASS90] is modified by coupling the non-attractive process to an attractive one.

*Proof of Proposition 2.7.* We take $\xi_t^-$ as a modified border contact process with parameters $(\lambda_c, \lambda_c + \varepsilon)$ and $\zeta_t^-$ a classical contact process with parameter $\lambda_c$. Denote by $\tilde{\mu}_t$ the law of $\xi_t^-$. Let $i, j$ be two fixed natural numbers and define

$$A_{i,j} = \{B \in \Sigma^{\circ} \; : \; |B \cap [-i, 0]| < j\}.$$

For every $k \in \mathbb{N}$, the random increment $\mathcal{R}\xi_{k+1} - \mathcal{R}\xi_k$ is stochastically dominated by a Poisson random variable of parameter $\lambda_c + \varepsilon$, thus,

$$\mathbb{E}[\mathcal{R}\xi_{k+1}^- - \mathcal{R}\xi_k^-]^+ \leq \lambda_c + \varepsilon. \tag{2.8}$$

When a translated version of $A_{i,j}$ occurs at time $s$, and at least $j$ death events and no birth events occur on $[\mathcal{R}\xi_t^- - i, \mathcal{R}\xi_t^-]$, then the right edge goes to the left by at least $i$ sites. Hence,

$$\mathbb{E}[\mathcal{R}\xi_{k+1}^- - \mathcal{R}\xi_k^-]^- \geq ip(j)\tilde{\mu}_k(A_{i,j}), \tag{2.9}$$

where $p(j)$ is a positive constant depending only on $j$. Multiplying (2.9), by $-1$ and summing with (2.8),

$$\mathbb{E}[\mathcal{R}\xi_{k+1}^- - \mathcal{R}\xi_k^-] \leq (\lambda_c + \varepsilon) - ip(j)\tilde{\mu}_k(A_{i,j}). \tag{2.10}$$

The processes $(\xi_t^-)_t$ and $(\zeta_t^-)_t$ can be coupled in an obvious way such that

$$\mathcal{R}\zeta_t^- \leq \mathcal{R}\xi_t^-.$$

Combining this with (2.10), we obtain

$$\begin{aligned}
\frac{1}{n}\sum_{k=0}^{n-1}[\mathbb{E}[\mathcal{R}\zeta_{k+1}^- - \mathcal{R}\zeta_k^-] &= \frac{\mathbb{E}[\mathcal{R}\zeta_n^-]}{n} \\
&\leq \frac{\mathbb{E}[\mathcal{R}\xi_n^-]}{n} \\
&= \frac{1}{n}\sum_{k=0}^{n-1}\mathbb{E}[\mathcal{R}\xi_{k+1}^- - \mathcal{R}\xi_k^-] \\
&\leq (\lambda_c + \varepsilon) - ip(j)\frac{1}{n}\sum_{k=0}^{n-1}\tilde{\mu}_k(A_{i,j}). \tag{2.11}
\end{aligned}$$

Rearranging the terms,

$$\frac{1}{n}\sum_{k=0}^{n-1}\tilde{\mu}_k(A_{i,j}) \leq \frac{\lambda_c+\varepsilon}{ip(j)} - \frac{1}{ip(j)n}\sum_{k=0}^{n-1}[\mathbb{E}[\mathcal{R}\zeta_{k+1}^- - \mathcal{R}\zeta_k^-] \qquad (2.12)$$

By Proposition 2.6,

$$\lim_n \frac{1}{n}\sum_{k=0}^{n}\tilde{\mu}_k(A_{i,j}) = \tilde{\mu}(A_{i,j}) \qquad (2.13)$$

and by [Lig05, Theorem 2.19],

$$\lim_n \frac{1}{n}\sum_{k=0}^{n}[\mathbb{E}[\mathcal{R}\zeta_{k+1}^- - \mathcal{R}\zeta_k^-] = \alpha(\lambda_c), \qquad (2.14)$$

and it is a classical result that $\alpha(\lambda_c) = 0$ (see [Gri81, Section 7]).

Therefore, taking the limit as $n$ goes to infinity in (2.12) and using (2.13) and (2.14),

$$\tilde{\mu}(A_{i,j}) \leq \frac{\lambda_c+\varepsilon}{ip(j)}. \qquad (2.15)$$

Taking the limit when $i$ goes to infinity in (2.15), we obtain

$$\tilde{\mu}\left(A \in \Sigma \ : \ |A| < j\right) = 0,$$

and thus $\tilde{\mu}(\Sigma^{\circ}) = 1$. $\qquad\qquad\square$

## 2.3    Non-survival at criticality

In this section, we prove Theorem 2.5. The proof follows the same steps of the one for the critical contact process in [BG90]. We give the main steps here for completeness.

*Proof of Theorem 2.5.* Let $\eta$ be a process with parameters $(\lambda_i, \lambda_e)$, $\lambda_i > \lambda_e$. We will prove that, if $\mathbb{P}_{\lambda_i,\lambda_e}(\eta_t^0 \neq \emptyset$ for all $t > 0) > 0$, then there is $\delta > 0$ such that, for every $\tilde{\lambda}_e, \tilde{\lambda}_i$ with $\max\{\lambda_i - \tilde{\lambda}_i, \lambda_e - \tilde{\lambda}_e\} < \delta$, $\mathbb{P}_{\tilde{\lambda}_i,\tilde{\lambda}_e}(\eta_t^0 \neq \emptyset$ for all $t > 0) > 0$. This is sufficient to prove the theorem.

**Lemma 2.16.** *Suppose that $\mathbb{P}_{\lambda_i,\lambda_e}(\eta_t^0 \neq \emptyset$ for all $t) > 0$. Then, for every $\varepsilon > 0$ there are $R, L \in \mathbb{N}$, $S > 0$ and $\delta > 0$ such that for every $\tilde{\lambda}_e, \tilde{\lambda}_i$ with $\max\{\lambda_i - \tilde{\lambda}_i, \lambda_e - \tilde{\lambda}_e\} < \delta$, with $\mathbb{P}_{\tilde{\lambda}_i,\tilde{\lambda}_e}$-probability at least $1 - \varepsilon$ there* *exists* $(x_0, t_0) \in [L, 2L] \times [S, 2S]$ *such that* $[-R, R] \times \{0\}$ *is connected inside* $[-L, 3L] \times [0, 2S]$ *to every point in* $[-R + x_0, R + x_0] \times \{t_0\}$.

We defer the proof of Lemma 2.16 to the end of this section.

The next step of the proof, as in [BG90], is to restart the process from the infected copy of $[-R, R]$. In contrast to what happens to the usual contact process, the contact process with modified border is not local, i.e., the probability of the spreading of the infection depends on the whole configuration. But, we can take $\delta$ sufficiently small such that $\tilde{\lambda}_i > \tilde{\lambda}_e$, so the process is still attractive and existence of infected particles outside $[-R + x_0, R + x_0] \times \{t_0\}$ at time $t_0$ can only increase the probability of infection. Thus, Lemma 2.16 can be iterated using the strong Markov property to obtain the following lemma.

**Lemma 2.17.** *Suppose that* $\lambda_i > \lambda_e$, $\mathbb{P}_{\lambda_i, \lambda_e}(\eta_t^0 \neq \emptyset \text{ for all } t) > 0$ *and* $R, L, S$ *as in Lemma 2.19. Then, for every* $\varepsilon > 0$, $k \in \mathbb{N}$, $x \in [-2L, 2L]$ *and* $t \in [0, 2S]$, *for every* $\tilde{\lambda}_e, \tilde{\lambda}_i$ *with* $\max\{\lambda_i - \tilde{\lambda}_i, \lambda_e - \tilde{\lambda}_e\} < \delta$, *with* $\mathbb{P}_{\tilde{\lambda}_i, \tilde{\lambda}_e}$*-probability at least* $(1 - \varepsilon)^{2k}$ *there exists* $(x_0, t_0) \in [L, 2L] \times [S, 2S]$ *such that* $[-x + R, x + R] \times \{t\}$ *is connected to every point in* $[-R + x_0, R + x_0] \times \{t_0\}$ *by paths lying entirely inside the region*

$$\bigcup_{j=0}^{k=1} [-3L + jL, 4L + jL] \times [2jS, 4S + 2jS].$$

*Proof.* We begin by proving that, with $\mathbb{P}_{\tilde{\lambda}_i, \tilde{\lambda}_e}$-probability at least $(1 - \varepsilon)^2$, for every $x \in [-2L, 2L]$ and $t \in [0, 2S]$ the interval $[x - R, x + R] \times \{t\}$ is connected to $[y - R, r + R] \times \{s\}$, for some $(y, s) \in [-L, 3L] \times [2S, 4S]$, and this connection is made through paths in the box $[-3L, 4L] \times [0, 4S]$. This step is necessary to ensure when applying the Lemma 2.16 repeatedly, small deviations of the "target" point do not accumulate and cause the desired infected disk to be centered outside the "target" region. Then, after the first application of Lemma 2.16 the direction of the path of infection is "corrected" if necessary.

There are five cases to consider, depending on the region of $[-2L, 2L] \times [0, 2S]$ in which the starting point $(x, t)$. The cases are pictured in Figure 2.2. In each case, the starting point is in the black region, and Lemma 2.16 (and symmetry) gives that with $\mathbb{P}_{\tilde{\lambda}_i, \tilde{\lambda}_e}$-probability at least $1 - \varepsilon$ every point in the black region is connected to some point in the gray region py paths lying entirely in the region enclosed by dashed lines. In cases $(iii)$ to $(v)$, the infected copy of $[-R, R]$ given by Lemma 2.16 may be centered outside $[-L, 3L] \times [2S, 4S]$. If this happens, we take the center of the earliest in time of those copies and apply case $(i)$ or $(ii)$ again as necessary.

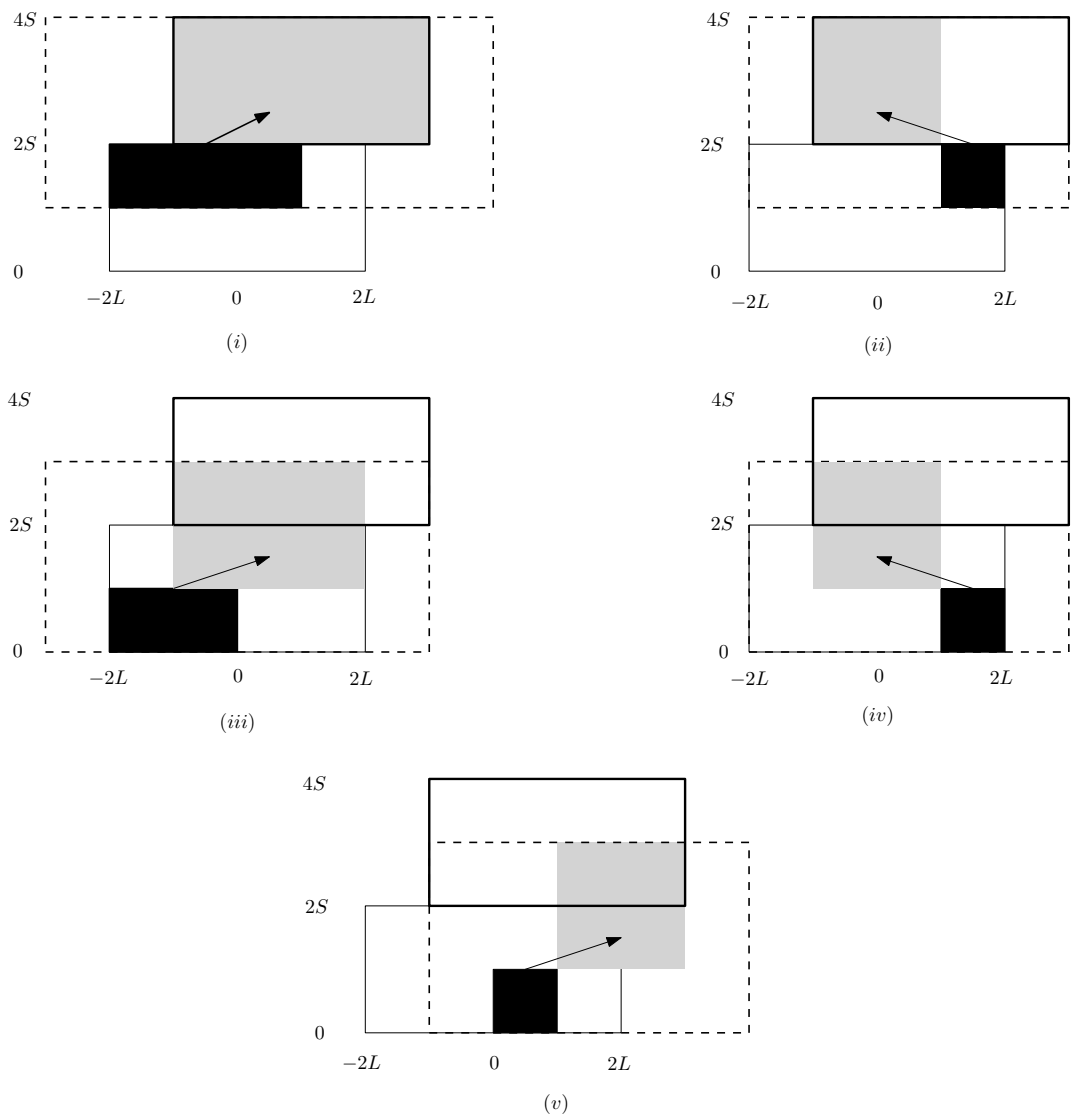

**Figure 2.2:** The cases that can occur when applying Lemma 2.16. At each case, the starting point is in the black area, and Lemma 2.16 ensures that the starting point is connected to some interval centered in the gray area.

In any case, we conclude that any infected copy of $[-R, R]$ centered $[-2L, 2L] \times [0, 2S]$ is connected to some infected copy of $[-R, R]$ centered in the box $[-L, 3L] \times [2S, 4S]$ with $\mathbb{P}_{\tilde{\lambda}_i, \tilde{\lambda}_e}$-probability at least $(1 - \varepsilon)^2$. Applying the strong Markov property $k$ times we conclude the proof of the lemma. $\square$

Fixing $k \in \mathbb{N}$, define the regions

$$V^{\pm} = \left\{ (x, t) \in \mathbb{Z} \times \mathbb{R} : 0 \leq t \leq (2k + 2)S, 5L \pm \frac{L}{2S}t \leq x \leq 5L \pm \frac{L}{5S}t \right\},$$

and, for $x \in [-2L, 2L]$ and $t \in [0, 2S]$ the event

$$G^{\pm}(x, t) = \left\{ \begin{array}{l} [-R + x, R + x] \times \{t\} \text{ is connected inside } V^{\pm} \text{ to every point in} \\ [-R + y, R + y] \times \{s\} \text{ for some } y \in [(\pm k - 2)L, (\pm k + 2)L] \\ \text{and } s \in [2KS, 2(k + 1)S]. \end{array} \right\}$$

Fix $\tilde{\varepsilon} > 0$ (to be fixed later), $k > 10$ and $\varepsilon$ such that $(1 - \varepsilon)^{2k} > 1 - \tilde{\varepsilon}$. By Lemma 2.17, for such $\varepsilon$ and $k$, there are $R$, $L$, $S, T, \delta$ such that, for all $(x, t) \in [-2L, 2L] \times [0, 2S]$, for every $(\tilde{\lambda}_i, \tilde{\lambda}_e)$ in the open ball of center $(\lambda_i, \lambda_e)$ and radius $\delta$, $\mathbb{P}_{\tilde{\lambda}_i, \tilde{\lambda}_e}(G^{\pm}(x, t)) > 1 - \tilde{\varepsilon}$. By a classic argument of 1-dependent site percolation (see [BG90, Lemma 21] for a detailed construction), if $\tilde{\varepsilon}$ is small enough, we have that

$$\mathbb{P}_{\tilde{\lambda}_i, \tilde{\lambda}_e}(\eta_t^0 \neq \emptyset \text{ for all } t > 0) > 0$$

and the proof is finished. $\square$

We prove now Lemma 2.16. The idea is to show existence of many infected paths to points in the top and the sides of a large space-time box, in a way that at least one of those points will result in an infected copy of the interval.

*Proof of Lemma 2.16.* We say that a point $(x, 0)$ is *connected to $+\infty$ by $(\lambda_i, \lambda_e)$-active open paths* if, for every $t > 0$ there is $y \in \mathbb{Z}$ such that $(x, 0)$ is connected to $(y, t)$ by $(\lambda_i, \lambda_e)$-active open paths for the initial configuration $\{x\}$.

By ergodicity of the measure $\mathbf{P}$, if the probability of survival is positive, there are an infinite number of $x \in \mathbb{Z}$ such that $(x, 0)$ connected to $+\infty$ by $(\lambda_i, \lambda_e)$-active open paths, $\mathbf{P}$-a.s.. Then, a disk centered $[-R, R]$ at the origin with diameter big enough will, with high probability, contain at least one of those points. By attractiveness, the probability that a contact process with initial condition $[-R, R]$ survives it grater than the probability that a contact process with initial condition $y \in [-R, R]$ survives. The, we can conclude that, for every $\varepsilon > 0$, there is $R > 0$ such that

$$\mathbb{P}(\eta_t^{[-R,R]} \neq \emptyset \text{ for all } t \geq 0) > 1 - \frac{1}{2}\varepsilon^4. \tag{2.18}$$

For $l \in \mathbb{N}$ and $t > 0$, define $N_T^+(l, t)$ (resp. $N_T^-(l, t)$) the number of points in $[-l, 0) \times t$ (resp. $[0, l] \times t$) which are joined to $[-R, R] \times 0$ inside $(-l, l) \times (0, t)$. Define

$$N_T(l, t) = N_T^+(l, t) + N_T(l, t)$$

to be the total number of points in the *top* of the space-time box $[-l, l] \times [0, t]$ which are joined to $[-R, R] \times 0$ inside the box. Analogously, we will establish a notation for the points infected in the *sides* of the space-time box. Choose $h \in (0, (1 + 2\lambda_e)^{-1})$, $l \in \mathbb{N}, t > 0$ and define $N_S^+(l, t)$ (resp. $N_S^-(l, t)$) as the maximal subset of $\{l\} \times [0, t]$ (resp $\{-l\}) \times [0, t]$ such that every two points are at $L^\infty$-distance of at least $h$; and $[-R, R] \times \{0\}$ is connected inside $[-l, l] \times [0, t]$ to every point of the set.

We then define

$$N_S(l, t) = N_S^-(l, t) + N_S^+(l, t)$$

and

$$N(l, t) = N_T(l, t) + N_S(l, t). \tag{2.19}$$

Take $a$ to be the minimum of

- the probability that the origin is connected inside $[-R, R] \times [0, h]$ to every point of $[-R, R] \times \{h\}$ by $(\lambda_i, \lambda_e)$-active open paths for the origin;

- the probability that the origin is connected inside $[0, 2R] \times [0, h]$ to every point of $[0, 2R] \times \{h\}$ by $(\lambda_i, \lambda_e)$-active open paths for the origin.

We then select $M$ big enough so that in $M$ independent Bernoulli trials of parameter $a$, the probability of at least one success is at least $1 - \varepsilon$; and then pick $N$ such that, in a collection of any given $N$ points in $\mathbb{Z}$, $M$ of them will be at distance at least $3R + 1$ apart.

Let $(\zeta_t)_t$ be a contact process with parameter $(\lambda_i, \lambda_i)$. By [DS00, Proposition 1], the probability of survival of $\zeta$ is greater than or equal to the probability of survival of $\zeta$. This observation, together with the attractiveness of the process and the arguments of the proof of [BG90, Lemma 7] show that there are $L \in \mathbb{N}$ and $T > 0$ such that

$$\mathbb{P}_{\lambda_i, \lambda_e}(N_T^{\pm}(L, T) \geq N) \geq 1 - \varepsilon \tag{2.20}$$

and

$$\mathbb{P}_{\lambda_i, \lambda_e}(N_S^{\pm}(L, T) \geq M) \geq 1 - \varepsilon. \tag{2.21}$$

By (2.21), with probability at least $1 - \varepsilon$ there are at least $M$ infected points at distance at least $h$ apart in $\{L\} \times [0, T]$. By the way $M$ was chosen, with probability at least $1 - \varepsilon$ one of those points (say, $(x, t)$) is connected by $(\lambda_i, \lambda_e)$ active paths in $[x+R, x+2R] \times [t, t+h]$ to every point in $[x+R, x+2R] \times \{t+h\}$. Define $\tau$ to be the first time such that $0$ is connected by active paths in $[-L, L+2R] \times [0, t]$ to every point in $[-L+2R, L+2R] \times \{t\}$. Then, we have that that the probability that $\tau$ is finite is at least $(1 - \varepsilon)^2$.

By the strong Markov property and Equation (2.20), conditioned on $\tau < +\infty$, with probability at least $1 - \varepsilon$, $[L, L+2R] \times \{\tau\}$ is connected to at least $N$ points in $[L+R, 2L+R] \times \{T+\tau\}$ by active paths lying inside $[L+R, 2L+R] \times [\tau, T+\tau]$.

By the way $N$ was chosen, at least $M$ of those $N$ points are at least a distance of $3R + 1$ apart. We divide those $M$ points in two sets:

- points $(z, \tau + T)$ with $z \leq 2L$. To each those points we associate the cylinder $[z - 2R, z] \times [T + \tau, T + \tau + h]$;

- points $(z, \tau + T)$ with $z > 2L$. To each those points we associate the cylinder $[z - R, z + R] \times [T + \tau, T + \tau + h]$.

All of those cylinders are disjoint. Using this fact and the definition of $M$, with probability at least $1 - \varepsilon$ one of those $M$ points, that we will call $x_0$, is joined inside its associated cylinder to every point on the top of the cylinder.

When this happens, there is a $(\lambda_i, \lambda_e)$-active path from $[-R, R] \times \{0\}$ to $[-R + x_0] \times \{2(T + h)\}$ lying entirely inside $[-L, 3L] \times [0, 2(T + h)]$ with probability at least $1 - 4\varepsilon$. Taking $S = T + h$, we conclude that there exists $(x_0, t_0) \in [L, 2L] \times [S, 2S]$ such that the origin is connected inside $[-L, 3L] \times [0, 2S]$ to every point in $[-R + x_0, R + x_0] \times \{t_0\}$ with $\mathbb{P}_{\lambda_i, \lambda_e}$-probability at least $1 - \varepsilon$. Since this event depends only on the configuration inside a finite space-time box, we can find $\delta > 0$ such that the open ball $B$ of radius $\delta$ and center $(\lambda_i, \lambda_e)$ is contained in the region $\{(\lambda_i, \lambda_e), : \lambda_i > \lambda_e \geq 0\}$ where the process is attractive, and the event still happens with $\mathbb{P}_{\tilde{\lambda}_i, \tilde{\lambda}_e}$-probability at least $1 - \varepsilon$ for every $(\tilde{\lambda}_i, \tilde{\lambda}_e) \in B$. $\qquad\square$

# Chapter 3

# Quasi-stationary distributions for subcritical population processes

Subcritical populational processes in general may not admit stationary distributions, but many of them have a *quasi-stationary distribution* (QSD), that is, a distribution that is invariant when conditioned on survival. In contrast to what generally happens for stationary distributions, QSDs may not be unique, even under irreducibility conditions. The characteristics of the process that may prevent it from possessing multiple QSDs are not entirely clear. A survey on the principal results about QSDs may be found on [MV12].

Formally, let $(\xi_t)_t$ be a Markov process on $\Lambda_0 := \Lambda \cup \{\emptyset\}$, with $\Lambda$ a countable set and $\emptyset$ an absorbing state that is reached a.s.. In all this chapter, if $\mu$ is a measure on $\Lambda$, we denote by $\mathbb{P}^\mu$ the law of the Markov process with initial condition sampled from $\mu$, and for $i \in \Lambda$, we write $\mathbb{P}^i$ instead of $\mathbb{P}^{\delta_i}$.

A measure $\nu$ on $\Lambda$ is a *quasi-stationary distribution* of the process $(\xi_t)_t$ if

$$\mathbb{P}^\nu(\xi_t \in \cdot \mid \xi_t \neq \emptyset) = \nu(\cdot),$$

for all $t > 0$.

The pioneer of the studies of quasi-stationary distributions is A. M. Yaglom. In [Yag47], he proved that, for a subcritical discrete-time branching process with finite variance $(X_n)_{n \geq 0}$, the limit

$$\lim_n \mathbb{P}^j(X_n \in \cdot \mid X_n \neq 0) := \nu(\cdot)$$

does not depend on $j \in \{1, 2, \dots\}$. Moreover, $\nu$ is a QSD on $\mathbb{N}$ for the branching process.

Seneta and Vere-Jones proved in [SVJ66] that there is a one-parameter family of QSDs for the subcritical branching process. Although presented for the discrete setting, the proof is the same for the continuous case [Cav78, Section 5]. Later, Cavender [Cav78] and Van Doorn [VD91] characterized the QSDs of general birth-and-death continuous time processes, proving that the family of [SVJ66] are the only QSDs for the subcritical branching process.

In contrast, uniqueness of QSD can be proved using some assumptions on the rate the process comes down from infinity (i.e., when the process, starting from an infinite configuration, reaches a finite configuration in a finite time) [BAJ24, CCL+09]. Although the subcritical contact process decays at the same rate as the branching process, it admits only one QSD, as shown in [AGR20]. Possible reasons for this difference may be because contact process, unlike the branching process, has a geometrical aspect. Thus, we may ask if it is the case that when adding geometrical structure to the branching process we recover the uniqueness property.

One way to do this is the *branching process with genealogy*, where we add information about the relationships between individuals. The state space is $\Sigma_0$, the set of finite rooted trees. Individuals are represented by vertices of the tree. The alive individuals are the leaves of the tree. Each individual, after a lifetime exponentially distributed, gives birth to a random (possibly zero) number of children and dies. The numbers of children of each individual are i.i.d. integrable random variables. Those children are attached to the tree by connecting them to their parent. To avoid unlimited growth of the tree, we delete the vertices and edges that are unnecessary to determine the genealogical relationship of the alive individuals. Thus, at each time $t$, we are able to We denote this process by $(\zeta_t)_t$. See Section 3.1 for a more precise description.

The projection $(\pi(\zeta_t))_t$ of the branching process with genealogy onto $\mathbb{N}_0$, with $\pi(\zeta_t)$ being the number of leaves of $\zeta_t$, is the usual branching process on $\mathbb{N}_0$ with the same offspring distribution of the branching process with genealogy. We note that if we start with a measure on $\Sigma_0$, evolve it by $t$ time units according with the dynamics of the branching process with genealogy and then project it onto $\mathbb{N}_0$, we obtain the same measure as if we start with a measure in $\Sigma_0$, project it onto $\mathbb{N}_0$ and then evolve it by $t$ time units using the dynamics of a branching process. In other words, the diagram of Figure 3.1 is commutative. This implies that, when the offspring distribution's mean is less than one, since the subcritical branching process has no non-trivial invariant measure, the same is true for the branching process with genealogy.

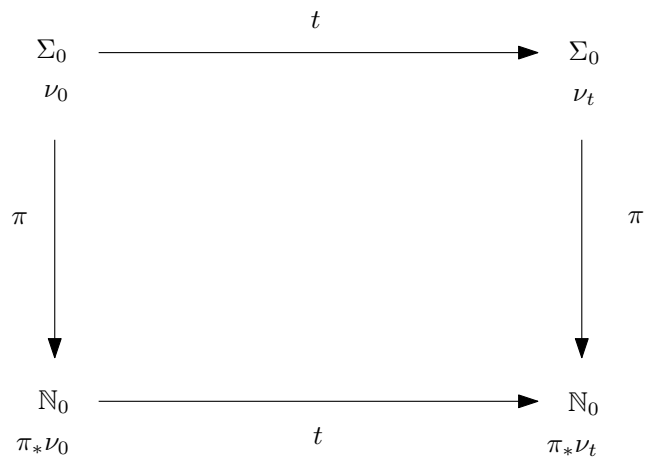

**Figure 3.1:** Relationships between measures

Define $\Sigma := \Sigma_0 \setminus \{\emptyset\}$, the set of non-empty finite rooted trees. Coming back to Figure 3.1 and noting that $\pi^{-1}\{0\} = \{\emptyset\}$ we see that we can project a measure on $\Sigma$ onto $\mathbb{N}$, then evolve for $t$ time units according with the dynamics of a branching process and condition on non-absorption, or we can start with a measure on $\Sigma$, evolve for $t$ time units using the dynamics of the branching process with genealogy, condition on non-absorption, and then project onto $\mathbb{N}$, and these two give the same result. Therefore, QSDs on $\Sigma$ give QSDs on $\mathbb{N}$ with the same absorption rate.

Assuming some condition on the moments of the offspring distribution, the branching process with genealogy admits at least one QSD. This is a consequence of the following proposition.

**Proposition 3.1.** *Let $(\xi_t)_t$ be an irreducible continuous-time Markov chain on a countable state-space $\Lambda_0 := \Lambda \cup \{\emptyset\}$ absorbed at $\emptyset$ a.s.. Suppose that there is a projection $\pi : \Lambda_0 \to \mathbb{N}_0$ such that $(\pi(\xi_t))_t$ is a branching process with offspring distribution $Z$ satisfying $\mathbb{E}e^{sZ} < +\infty$ for some $s > 0$, and $\mathbb{P}(Z = 0) + \mathbb{P}(Z = 1) < 1$. Suppose also that $\pi^{-1}\{0\} = \{\emptyset\}$ and that there is $\mathbf{1} \in \Lambda$ such that $\pi^{-1}\{1\} = \{\mathbf{1}\}$. Then there is a QSD $\nu$ of $(\xi_t)_t$ on $\Lambda$ such that,* *for every $\xi_0, \xi \in \Sigma$,*

$$\lim_{t \to +\infty} \mathbb{P}^{\xi_0}(\xi_t = \xi \mid \xi_t \neq \emptyset) = \nu(\xi). \tag{3.2}$$

A QSD that satisfies the limit (3.2) is called a *Yaglom limit*.

Contrary to what happens with the projection onto $\mathbb{N}$, we prove that the Yaglom limit is the unique QSD for the branching process with genealogy. This is one of the two main new contributions of this chapter.

**Theorem 3.3.** *Let $(\zeta_t)_t$ be a branching process with genealogy with offspring distribution $Z$ satisfying $\mathbb{E}[Z] < 1$, $\mathbb{P}(Z = 0) + \mathbb{P}(Z = 1) < 1$ and $\mathbb{E}^{sZ} < +\infty$ for some $s > 0$. Assume also that $\Sigma$ is an irreducible class of $(\zeta_t)_t$. Then $(\zeta_t)_t$ has a unique quasi-stationary distribution $\nu$ on $\Sigma$.*

The main idea behind the proof of Theorem 3.3 is to show that the geometrical aspects of the process prevent the existence of more than one QSD. To be more precise, we prove that, conditioned on non-absorption at time $t$, the number of branchings goes to infinity. Once this is done, we investigate the descendants of the first leave at time $t = 0$ that have descendants alive at time $t$. Either this set of descendants is the entire tree at time $t$ (and thus is distributed approximately as the Yaglom limit), or the tree at time $t$ has to be very "large", contradicting the fact that the distribution is quasi-stationary.

Another way to include a geometrical aspect to the branching process is the *branching random walk*. A branching random walk is a stochastic process that includes characteristics of both branching processes and random walks. Define $\Delta^d := \{\eta \in \mathbb{N}_0^{\mathbb{Z}^d} : \eta(x) \neq 0 \text{ only for a finite number of } x \in \mathbb{Z}^d\}$. We describe the branching random walk as a process on $\Delta_0^d$ in the following way. A finite number of particles is distributed on the lattice $\mathbb{Z}^d$. Let $\lambda > 0$ be a parameter. When there is a particle $X$ at site $x$, this particle chooses, at rate $\lambda$, one of its neighbors $y$ uniformly at random and give birth to a new particle $Y$ in $y$. The particle $Y$ is called a *child* of particle $X$. A particle in $\mathbb{Z}^d$ waits an exponential time with rate 1 independently of all toher particles and then dies. We define $\eta_t^A$ as the configuration at time $t$ starting with the configuration $A$ at time 0. The superscript may be omitted when the initial configuration is clear from context. There is a huge literature about branching random walks. Some introductory references may be found in [Sch99, Shi15].

The process as described above does not admit a QSD. Conditioned on survival on time $t$, the configuration $\eta_t^A$ is, in general, very far from the location of the initial configuration $A$. Therefore, we work with the branching random walk *modulo translations*. Take the equivalence relationship on $\mathbb{N}_0^{\mathbb{Z}^d}$ in which $\eta \sim \eta'$ if $\eta'$ is a translation of $\eta$. Denote by $\langle \eta \rangle$ the equivalence class of $\eta$. The branching random walk modulo translations is the process $(\langle \eta_t \rangle)_t$ on $\Delta_\sim^d := \Delta^d / \sim$.

We note that the projection of this process on $\mathbb{N}_0$ is a branching process with offspring distribution supported on $\{0, 2\}$ and the critical parameter is $\lambda = 1$. Thus, by Proposition 3.1, there it exists at least one QSD $\nu'$ in $\Delta_\sim^d$ for the branching random walk when $\lambda < 1$.

The branching random walk can be seen as a contact process with multiplicities, i.e., a contact process where it is possible to infect again a site which is already infected. Similarly to what happens to the branching process with genealogy, its geometrical aspects impede the existence of multiple QSDs. This is our second main contribution in this chapter.

**Theorem 3.4.** *If $\lambda < 1$, there is only one QSD, $\nu'$, in $\Delta_\sim^d$ for the branching random walk modulo translations.*

The proof of Theorem 3.4 relies on techniques similar to those of Theorem 3.3. We look for an individual at time $s = 0$ that has alive descendants at time $s = t$. Either this set of descendants are the only individuals alive at time $t$, and has a distribution close to the Yaglom limit, or there are another individual at time $s = 0$ with alive descendants at time $t$. But the sets of descendants of those two individuals are in general very distant from each other at time $t$, and thus the population at time $t$ would be spread over a very large set. Again, this would be incompatible with the concept of a QSD.

To prove that conditioned on survival, the number of births in the branching process with genealogy and the branching random walk will be very large, and thus conclude the impossibility of existence of multiple QSDs, we will need the following proposition.

**Proposition 3.5.** *Let $(X_t)_t$ be a subcritical branching process with offspring distribution satisfying the same hypothesis as in Proposition 3.1. Define, for $k \in \mathbb{N}$ and $t > 0$,*

$$G_t^X(k) = \left\{ \begin{array}{l} \text{There are } 0 < t_1 < \cdots < t_{2k} < t : \\ X_{t_1} = 1, X_{t_2} = 2, \ldots, X_{t_{2k-1}} = 1, X_{t_{2k}} = 2 \end{array} \right\}.$$

*Then, for every $k \in \mathbb{N}$,*

$$\lim_{t \to +\infty} \mathbb{P}^1(G_t^X(k) \mid X_t \neq 0) = 1.$$

The use of state 2 in the above proposition is not essential. In fact, 2 could be replaced with any integer greater than 1. The important consequence is that Proposition 3.5 implies that the number of branching events on a subcritical branching process will exceed any predetermined number, when conditioned on survival. The proof of Proposition 3.5 will be presented at Section 3.5.

The remaining of this chapter is organized as follows. In Section 3.1 we give a formal description of the branching process with genealogy. Theorems 3.3 and 3.4 are proved using the arguments outlined above in Sections 3.2 and 3.3 respectively. Proposition 3.1 is proved using classical results about $\alpha$-positiveness of submarkovian kernels in Section 3.4. Finally, in Section 3.5 we prove that conditioned on survival, there are many branching events by approximating the branching process by a process that is conditioned on not dying for a long time.

## 3.1   The branching process with genealogy

In this section, we formally define the branching process with genealogy. A classical construction of the branching process with family trees can be found e.g. on [Har02, Chapter VI].

Let $\Sigma$ be the set of non-empty finite rooted trees, and $\Sigma_0 = \Sigma \cup \{\emptyset\}$. We consider two trees to be equivalent if there is a graph isomorphism that maps one onto the other preserving the root. So, $\Sigma$ is a countable space.

A vertex $x$ is a *descendant* of a vertex $y$ or, equivalently, $y$ is an *ancestor* of $x$ if there is an oriented path from $y$ to $x$, the positive direction being away from the root. By convention, we say that every vertex $x$ is an ancestor and a descendant of itself.

The *diameter* of a tree $\zeta$ is the length of the greatest (non-oriented) path of the tree, and is denoted by $\text{diam}(\zeta)$. We denote by $\mathbf{1}$ as the tree with only one vertex.

We construct the branching process with genealogy $(\zeta_t)_t$ as a Markov process

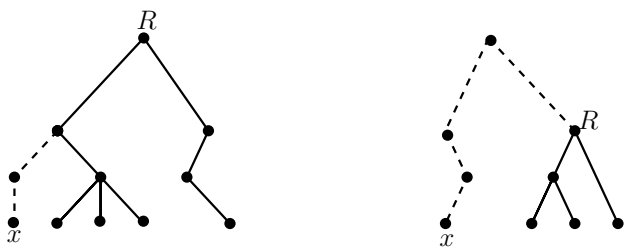

**Figure 3.2:** Two examples of pruning. The dashed piece of the tree is deleted
when vertex $x$ dies.

on $\Sigma_0$. The dynamics of the process is the following. Vertices can be either *alive*
or *dead*. All the leaves are alive and the other vertices are dead. Each alive
vertex $x$ carries an exponential clock that rings with rate 1 independently of
each other. When the clock of vertex $x$ rings, we sample a independent copy
$Z^x$ of a non-negative integer-valued random variable $Z$, called *the offspring
distribution*. Then, a number $Z^x$ of vertices is added to the tree, each one of
them connected to vertex $x$. When this happens, we say that the vertex $x$
*gave birth*, or that a *branching event* occurred, and the state of $x$ is changed
from alive to dead.

Every time a vertex dies, we *prune* the tree. That means we remove from
the tree all the vertices that are unnecessary to determine the affinity between
every pair of alive vertices. To be more precise, suppose that the Poisson clock
of vertex $x$ rings at time $t$. If $x$ gives birth to at least one child, we do nothing.
If $x$ gives birth to zero children and die, we prune the tree. Call $R$ the most
recent common ancestor to all alive leaves at time $t$ (except for $x$). That is,
all alive vertices at time $t$ descend from $R$ and $R$ is the farthest vertex from
the root with this property. There are two possible cases:

(i) $R$ is the root. In this case, we delete all the ancestors of $x$ which do not
    have other alive descendants at time $t$.

(ii) $R$ is not the root. In this case, we delete all the vertices which are not
     descendants of $R$, and declare $R$ as the new root.

Both cases are pictured in Figure 3.2.

If $\nu$ is a measure on $\Sigma$, we denote by $\mathbb{P}^\nu$ the law of this process with initial
condition distributed as $\nu$. We write $\mathbb{P}^A$ for $\mathbb{P}^{\delta_A}$.

## 3.2 The QSD for the branching process with genealogy

In this section we prove Theorem 3.3.

*Proof of Theorem 3.3.* Fix $t > 0$. We define the event

$$\mathcal{O}_t = \{\text{There is only one leaf at } \zeta_0 \text{ with alive descendants at time } t\}$$

Define for every leaf $Y$ of $\zeta_0$, $\mathcal{D}_t^Y$ as the subtree of $\zeta_t$ induced by the descendants of $Y$. That is, $\mathcal{D}_t^Y$ is the graph whose vertices are the descendants of $Y$ that are in $\zeta_t$ and the edges of $\zeta_t$ connecting those vertices.

Let $\nu^*$ be a QSD in $\Sigma$ of $(\zeta_t)_t$. For every $\zeta \in \Sigma$,

$$\begin{aligned}
\nu^*(\zeta) &= \mathbb{P}^{\nu^*}(\zeta_t = \zeta \mid \zeta_t \neq \emptyset) \\
&= \mathbb{P}^{\nu^*}(\zeta_t = \zeta \mid \mathcal{O}_t)\mathbb{P}^{\nu^*}(\mathcal{O}_t \mid \zeta_t \neq \emptyset) \\
&\quad + \mathbb{P}^{\nu^*}(\zeta_t = \zeta \mid \mathcal{O}_t^c \cap \{\zeta_t \neq \emptyset\})\mathbb{P}^{\nu^*}(\mathcal{O}_t^c \mid \zeta_t \neq \emptyset). \quad (3.6)
\end{aligned}$$

Conditioning on $\mathcal{O}_t$, the distribution of $\zeta_t$ is the distribution of a branching process with genealogy with initial configuration being a tree with a single leaf (which we denote by $\mathbf{1}$) conditioned on $\{\zeta_t \neq \emptyset\}$. Then,

$$\mathbb{P}^{\nu^*}(\zeta_t = \zeta \mid \mathcal{O}_t) = \mathbb{P}^{\mathbf{1}}(\zeta_t = \zeta \mid \zeta_t \neq \emptyset) \to \nu(\zeta), \quad (3.7)$$

where the last limit comes from Proposition 3.1.

On the event $\{\zeta_t \neq \emptyset\}$, there is at least one leaf at time 0 with alive descendants at time $t$. Choose uniformly at random a leaf $\mathcal{X}$ at time 0 which has an alive descendant at time $t$. Conditioned on $\{\zeta_t \neq \emptyset\}$, $\mathcal{D}_t^{\mathcal{X}}$ has the distribution of a branching process with genealogy starting from $\mathbf{1}$ conditioned on survival up to time $t$. Thus, remembering that $\pi(\mathcal{D}_t^{\mathcal{X}})$ is the number of leaves of $D_t^{\mathcal{X}}$, conditioned on $\{\zeta_t^\nu \neq \emptyset\}$, the distribution of the projection $\pi(\mathcal{D}_t^{\mathcal{X}})_t$ is the one of a branching process with initial state 1 conditioned on survival up to time $t$. Moreover, if $G_t^{\pi(\mathcal{D}_t^{\mathcal{X}})}(k) \cap \mathcal{O}_t^c \cap \{\zeta_t \neq \emptyset\}$ occurs, then $\zeta_t$ has diameter at least $k$. Indeed, when there are $t', t''$ such that $\pi(\mathcal{D}_{t'}^{\mathcal{X}}) = 1$ and $\pi(\mathcal{D}_{t''}^{\mathcal{X}}) = 2$ there is some $t \in (t', t'')$ when an edge is added to $\zeta_t$. Then, when $k$ of those transitions occur, at least $k$ edges have been added. On $\mathcal{O}_t^c \cap \{\zeta_t \neq \emptyset\}$, there are individuals alive at time $t$ which are not descendants of $\mathcal{X}$. By the rules to prune the tree, this implies that $\mathcal{X}$ and those $k$ added edges are still on $\zeta_t$ (they are necessary to establish the genealogical relationship between the individuals that are descendants of $\mathcal{X}$ and those who are not).

By Proposition 3.5, for every $\varepsilon > 0$ there is $t$ big enough such that $\mathbb{P}^1(G_t^{\pi(\mathcal{D}_t^{\mathcal{X}})}(k)) > 1 - \varepsilon$. Using this fact and the comments of the previous paragraph, we conclude that, for every $t$ big enough, $\mathbb{P}(\mathrm{diam}(\zeta_t) \geq k; \mathcal{O}_t^c \mid \zeta_t \neq \emptyset) > 1 - \varepsilon$. Then, for $t$ big enough we have that

$$\begin{aligned}
\mathbb{P}^{\nu^*}(\mathcal{O}_t^c \mid \zeta_t \neq \emptyset) &= \mathbb{P}^{\nu^*}(\{\mathrm{diam}(\zeta_t) \geq k\} \cap \mathcal{O}_t^c \mid \zeta_t \neq \emptyset) \\
&\quad + \mathbb{P}^{\nu^*}(\{\mathrm{diam}(\zeta_t) < k\} \cap \mathcal{O}_t^c \mid \zeta_t \neq \emptyset) \\
&\leq \mathbb{P}^{\nu^*}(\mathrm{diam}(\zeta_t) \geq k \mid \zeta_t \neq \emptyset) + \varepsilon \\
&= \nu^*\{\zeta : \mathrm{diam}\,\zeta \geq k\} + \varepsilon.
\end{aligned}$$

Letting $k \to +\infty$ and $\varepsilon \to 0$, we get

$$\lim_{t \to +\infty} \mathbb{P}^{\nu^*}(\mathcal{O}_t^c \mid \zeta_t \neq \emptyset) = 0. \tag{3.8}$$

Taking the limit in $t$ in (3.6) and using (3.7) and (3.8), we have that $\nu(\zeta) = \nu^*(\zeta)$. $\qquad\square$

## 3.3   Uniqueness of QSD for branching random walks

In this section, we prove uniqueness of the quasi-stationary distribution for the branching random walk. To do this, we argue that, if there are at time 0 two individuals with alive descendants at time $t$, the set of descendants of those two individuals are likely very far from each other, contradicting the fact that the diameter of a QSD distributed configuration cannot grow without bounds.

*Proof of Theorem 3.4.* We remember that $\Delta^d$ is the subset of $\mathbb{N}_0^{\mathbb{Z}^d}$ with only a finite number of sites with more than zero particles, $\sim$ the relation of equivalence on $\Delta^d$ with $\eta \sim \eta'$ if $\eta$ is a translation of $\eta'$ and $\Delta_\sim^d$ is the quotient of $\Delta^d$ by $\sim$. Usually, notations with $\sim$ refers to objects and events in the quotient space $\Delta_\sim^d$.

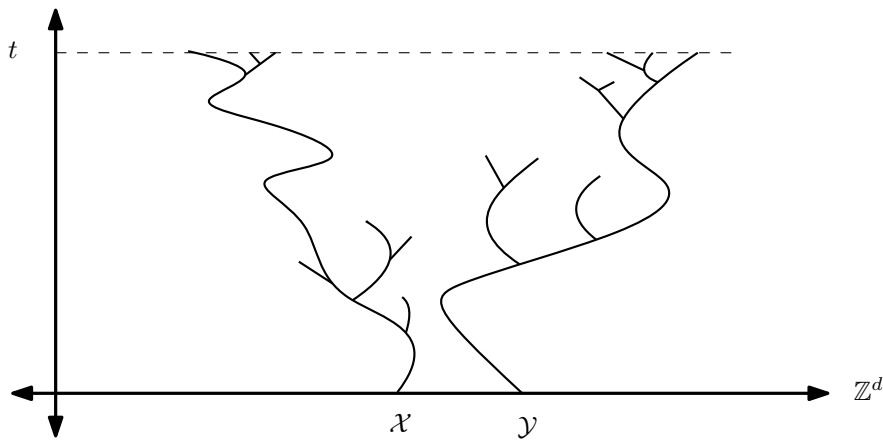

**Figure 3.3:** If there are two particles at time 0 with alive descendants at time $t$, there are particles very far from each other at time $t$

For a configuration $\eta \in \Delta^d$, define $\mathrm{diam}(\eta) = \mathrm{diam}(\{x \in \mathbb{Z}^d : \eta(x) \geq 1\})$. This is constant in each equivalence class, then we can define $\mathrm{diam}(\langle \eta \rangle) := \mathrm{diam}(\eta)$.

A particle $Y$ is said to be a *descendant* of a particle $X$ if there are particles $X = X_0, X_1, \dots, X_k = Y$ such that $X_j$ is a child of $X_{j-1}$ for $1 \leq j \leq k$. By convention, we always say that a particle is a descendant of itself. Denote, for each $t > 0$, as $\tilde{\mathcal{O}}_t$ the event in which there is exactly one individual at time 0 with alive descendants at time $t$. Define, for a particle $X$, $\mathcal{D}_t^X$ to be the restriction of the process $\eta$ to the descendants of particle $X$ at time $t$; that is, $\mathcal{D}_t^X(x)$ is the number of particles at time $t$ on site $x$ that are descendants of $X$.

We also define **0** to be the configurations with no particles.

Fix $t > 0$. On the event $\{\eta_t \neq \mathbf{0}\}$, we associate to each particle $\mathcal{X}$ alive at $t = 0$ that has alive descendants at time $t$ a particle $(S_{\mathcal{X}}(s))_{s \leq t}$, in the following way. At each time $s \leq t$ we call a particle *walker of $\mathcal{X}$ in the process* $(\eta_s)_{s \in [0,t]}$. At time $t = 0$, take $S_{\mathcal{X}}(0) = \mathcal{X}$. Suppose $(S_{\mathcal{X}}(s))_{s \in [0,w]}$ is defined, for some $w \leq t$, and define $\mathcal{Y} := S_{\mathcal{X}}(w)$. Let $\tau := \inf\{s \in (w, t]; \mathcal{Y} \text{ gives birth at time } s\}$. Call $\mathcal{Z}$ the particle born from $\mathcal{Y}$ at time $\tau$. We have the following cases:

- if $\mathcal{Z}$ has alive descendants at time $t$, then $S_{\mathcal{X}}(s) = \mathcal{Y}$ for $w < s < \tau$ and $S_{\mathcal{X}}(\tau) = \mathcal{Z}$. In this case, if $\mathcal{Y}$ is at site $x$ and $\mathcal{Z}$ is at site $x \pm e_i$, for some $1 \leq i \leq d$, we say that the walker *jumped* in direction $\pm e_i$, or just *jumped*, if it is not necessary to make the direction explicit;

- otherwise, we define $S_{\mathcal{X}}(s) = \mathcal{Y}$, for $w \leq s \leq \tau$.

In words, every time the walker $\mathcal{Y}$ gives birth to a new particle $\mathcal{Z}$, this newborn particle $\mathcal{Z}$ is called the walker if, and only if, $\mathcal{Z}$ has alive descendants at time $t$. Otherwise, $\mathcal{Y}$ continues to be the walker.

The following lemma ensures that the number of jumps of the walker process goes to infinity in probability, conditioned on survival.

**Lemma 3.9.** *For every $M \in \mathbb{N}$, and every $\eta_0 \in \Delta^d$ with only one particle,*

$\lim_{t \to +\infty} \mathbb{P}^{\eta_0}(\text{The walker of } X \text{ jumps at least } M \text{ times before time } t \,|\, \eta_t \neq \mathbf{0}) = 1,$

*where $X$ is the only particle alive at time $t = 0$.*

Conditioning on occurrence of $\tilde{\mathcal{O}}_t^c$ and survival up to time $t$, choose uniformly at $\mathcal{X}$ and $\mathcal{Y}$ two distinct individuals with alive descendants at time $t$. We have that $\{\eta_t \neq \mathbf{0}\} \cap \tilde{\mathcal{O}}_t^c \subseteq \{\text{diam}\langle\eta_t\rangle \geq |S_{\mathcal{X}}(t) - S_{\mathcal{Y}}(t)|\}$, for every possible choice of $\mathcal{X}$ and $\mathcal{Y}$.

Conditioned on $\{\eta_t \neq \mathbf{0}\} \cap \tilde{\mathcal{O}}_t^c$ and the jump times of $(S_{\mathcal{X}})_{s \in [0,t]}$, the distribution of the direction of each jump of the walker of $\mathcal{X}$ is independent of the direction of the other jumps and uniform in $\{\pm e_1, \pm e_2, \ldots, \pm e_d\}$. Therefore, conditioned on survival and on jump times, the jumps of the walker of $\mathcal{X}$ are those of a simple symmetric random walk. As the process at time $t$ conditioned of survival is the superposition of the independent process of the descendants of the particles alive at $t = 0$, the same is valid for the jumps of the walker of $\mathcal{Y}$, which are also conditionally independent of the jumps of $\mathcal{S}_{\mathcal{X}}$.

By Proposition 3.1, there is a QSD, $\nu'$, on $\Delta_{\sim}^d$ for the BWR modulo translation that satisfies the Yaglom limit. Let $\tilde{\nu}$ be an arbitrary QSD on $\Delta_{\sim}^d$ for the branching random walk modulo translation.

Let $\varepsilon > 0$. Choose $M \in \mathbb{N}$ satisfying the following conditions:

- for every $m \geq M$, the probability that the distance at time $m$ between two discrete-time independent simple symmetric random walks starting from arbitrary positions is less than $M^{1/4}$ is at most $\varepsilon$;

- $\tilde{\nu}\{\text{diam}\,\eta \geq M^{1/4}\} < \varepsilon$.

Denote by $\mathbf{1}$ the configuration on $\Delta^d_{\sim}$ with only one particle. Conditioned on $\tilde{O}_t$, $\langle \eta_t \rangle$ is distributed as $\langle \eta_t^{\mathbf{1}} \rangle$.

On the event $\tilde{\mathcal{O}}^c_t \cap \{\eta_t \neq \mathbf{0}\}$, the projection on on $\Delta^d_{\sim}$ of each one of the processes $(\mathcal{D}^{\mathcal{X}}_s)_{s \in [0,t]}$ and $(\mathcal{D}^{\mathcal{Y}}_s)_{s \in [0,t]}$ is distributed as a branching random walk modulo translation with initial configuration $\mathbf{1}$ conditioned on survival up to time $t$. By Lemma 3.9, for all $t_0$ big enough, conditioned on survival up to time $t_0$ and on $\tilde{\mathcal{O}}_{t_0}$, with probability at least $1 - \varepsilon$, the walkers of $\mathcal{X}$ and $\mathcal{Y}$ both make at least $M$ jumps before time $t_0$ (call this event $H^M_{t_0}$). Then, for $t > t_0$,

$$
\begin{aligned}
\mathbb{P}^{\tilde{\nu}}(\tilde{\mathcal{O}}^c_t \mid \eta_t \neq \mathbf{0}) &\leq \mathbb{P}^{\tilde{\nu}}(H^M_{t_0} \mid \tilde{\mathcal{O}}^c_t \cap \{\eta_t \neq \mathbf{0}\}) + \varepsilon \\
&\leq \mathbb{P}^{\tilde{\nu}}(|S_{\mathcal{X}}(t) - S_{\mathcal{Y}}(t)| \geq M^{1/4} \mid \tilde{\mathcal{O}}^c_t \cap \{\eta_t \neq \mathbf{0}\}) + 2\varepsilon \\
&\leq \mathbb{P}^{\tilde{\nu}}(\mathrm{diam}(\langle \eta_t \rangle) \geq M^{1/4} \mid \tilde{\mathcal{O}}^c_t \cap \{\eta_t \neq \mathbf{0}\}) + 2\varepsilon \\
&= \tilde{\nu}\{\mathrm{diam}(\eta) \geq M^{1/4}\} + 2\varepsilon \\
&\leq 3\varepsilon.
\end{aligned}
$$

As $\varepsilon$ is arbitrary,

$$
\lim_{t \to +\infty} \mathbb{P}^{\tilde{\nu}}(\tilde{\mathcal{O}}^c_t \mid \eta_t \neq \mathbf{0}) = 0. \tag{3.10}
$$

Then, for every $\eta \in \Delta^d_{\sim}$,

$$
\begin{aligned}
\tilde{\nu}(\eta) &= \mathbb{P}^{\tilde{\nu}}(\langle \eta_t \rangle = \eta \mid \tilde{\mathcal{O}}_t)\mathbb{P}^{\tilde{\nu}}(\tilde{\mathcal{O}}_t \mid \eta_t \neq \mathbf{0}) \\
&\quad + \mathbb{P}^{\tilde{\nu}}(\langle \eta_t \rangle = \eta \mid \tilde{\mathcal{O}}^c_t \cap \{\eta_t \neq \mathbf{0}\})\mathbb{P}^{\tilde{\nu}}(\tilde{\mathcal{O}}^c_t \mid \eta_t \neq \mathbf{0}) \\
&= \mathbb{P}^{\mathbf{1}}(\langle \eta_t \rangle = \eta \mid \eta_t \neq \mathbf{0})\mathbb{P}^{\tilde{\nu}}(\tilde{\mathcal{O}}_t \mid \eta_t \neq \mathbf{0}) \\
&\quad + \mathbb{P}^{\tilde{\nu}}(\langle \eta_t \rangle = \eta \mid \tilde{\mathcal{O}}^c_t \cap \{\eta_t \neq \mathbf{0}\})\mathbb{P}^{\tilde{\nu}}(\tilde{\mathcal{O}}^c_t \mid \eta_t \neq \mathbf{0}).
\end{aligned}
$$

Taking the limit in $t$, using (3.10) and the fact that $\nu'$ is a Yaglom limit, we obtain $\tilde{\nu} = \nu'$. $\qquad \square$

*Proof of Lemma 3.9.* Let $\varepsilon > 0$. Choose $M_1$ such that, for every $n \geq M_1$ the probability of at least $n$ successes in $3n$ Bernoulli trials with parameter $1/2$ exceeds $1 - \varepsilon/2$, and let $M_2 = \max\{3M, 3M_1\}$.

By Proposition 3.5, there is $t_0 > 0$ such that conditioned on survival at time $t_0$, the process $(D^X_t)_t$ goes from a state with only one particle to a state with two particles at least $M_2$ times.

The above choice of $t_0$ implies that, conditioned on survival up to time $t_0$, the probability that there are at least $M_2$ times when the walker gives birth to a particle exceeds $1 - \varepsilon/2$. That is, conditioned on $\{\eta^1 \neq \mathbf{0}\}$, the event $\{$There are $0 < t_1 < \cdots < t_{M_2} < t : S_{\mathcal{X}}(t_i^-)$ gives birth at time $t_i, 1 \leq i \leq M_2\}$ is at least $1 - \varepsilon/2$. Indeed, every time the configuration has only one particle this particle must be the walker; and if at a given time the configuration has more than one particle, this means that there must have been a birth involving the walker at some earlier time.

Suppose the walker $\mathcal{Z}$ of $X$ gives birth at a time $t_1 < t_0$ to a particle $\mathcal{W}$. By the branching property, the processes $(\mathcal{D}_s^{\mathcal{Z}})_{s \geq t_1}$ and $(D_s^{\mathcal{W}})_{s \geq t_1}$ (i.e., the processes of the descendants of $\mathcal{Z}$ and $\mathcal{W}$), conditioned on $\{\eta_{t_0} \neq \emptyset\}$ are equally distributed branching random walks. As $\mathcal{Z}$ is the walker, we know that at least one of those processes is alive at time $t_0$, and by symmetry, conditioning on $\sigma(\mathcal{D}_s^{\mathcal{X}})_{s \leq t_1}$ and that the walker at time $t_1$ is $\mathcal{Z}$, the probability that $\mathcal{D}_t^{\mathcal{W}} \neq \emptyset$ is at least $1/2$. Thus, each time a walker gives birth, the conditional probability that the walker makes a jump is at least $1/2$.

By the way $t_0$ was chosen, there is more than $M_2$ times the walker gives birth until time $t_0$ with probability at least $1 - \varepsilon/2$ (conditioned on survival). Also conditioning on survival, by the arguments of the previous paragraph and the choice of $M_1$, there is a probability at least $1 - \varepsilon/2$ that at least $M$ of those birth events results on a jump of the walker, and the lemma is proved. $\qquad \square$

## 3.4 $\alpha$-positiveness and the Yaglom limit

In this section we prove Proposition 3.1. In all this section, $(\xi_t)_t$ is an irreducible Markov chain on $\Lambda_0 := \Lambda \cup \{\emptyset\}$ with $\emptyset$ an an absorbing state reached a.s. and $\Lambda$ a countable set.

Let $P_t$ be the sub-Markovian kernel associated with the restriction of the process $\zeta$ to $\Lambda$. That is, for $i, j \in \Lambda$, $P_t(i,j) = \mathbb{P}(\xi_t^i = j)$. Then, by [Kin63, Theorem 1], for every $i, j \in \Lambda$ the limit

$$\alpha := \lim_{t \to +\infty} -\frac{\log P_t(i,j)}{t}$$

exists and does not depend on the choice of $i$ and $j$. Furthermore, $\alpha > 0$.

We say that the semi-group $(P_t)_t$ (or, equivalently, the process $(\xi_t)_t$) is $\alpha$-*recurrent* if

$$\int_0^{+\infty} e^{\alpha t} P_t(i,i) \, \mathrm{d}\, t > 0$$

for some (and thus, for all) $i \in \Lambda$, and $\alpha$-*positive* if

$$\limsup_{t \to +\infty} e^{\alpha t} P_t(i,i) > 0,$$

holds for some (equivalently, for all) $i \in \Lambda$.

We need the following result.

**Theorem 3.11** ([Kin63, Theorem 4]). *If $(P_t)_t$ is $\alpha$-recurrent, then there are positive vectors $\nu = (\nu_i)_{i \in \Lambda}$ and $h = (h_i)_{i \in \Lambda}$, unique up to a multiplicative constant, such that, for all $t \geq 0$,*

$$P_t h = e^{-\alpha t} h; \quad \nu P_t = e^{-\alpha t} \nu. \tag{3.12}$$

*The kernel $(P_t)_{t \geq 0}$ is $\alpha$-positive if, and only if, $\nu h < +\infty$. In this case, for all $i, j$ in $\Lambda$,*

$$\lim_{t \to +\infty} e^{\alpha t} P_t(i,j) = \frac{h_i \nu_j}{\nu h}. \tag{3.13}$$

If, in addition, the left-eigenvector $\nu$ is summable, existence of the Yaglom limit follow from classical arguments [SVJ66]. To prove summability of $\nu$, we make use of a discrete process. Let $(X_n)_n$ be a discrete-time, irreducible and aperiodic Markov chain on $\Lambda_0$, absorbed a.s. at $\emptyset$, with transition matrix $P(\cdot, \cdot)$. Analogously to the continuous-time case have that, for $i, j \in \Lambda$, $\lim_n (P_n(i,j))^{-1/n} = R > 1$, and we say that $(X_n)_n$ is $R$-positive if $\limsup R^n P_n(i,j) > 0$.

Define, for $i \in \Lambda$, $\tau^i := \inf\{n \in \mathbb{N} : X_n^i = \emptyset\}$, with $X_k^i$ being the chain $(X_n)_n$ starting from state $i$.

We make use of the following theorem.

**Theorem 3.14** ([FKM96, Theorem 1]). *Suppose there exists a subset $\Lambda' \subseteq \Lambda$, a configuration $A' \in \Lambda'$, $\rho < R^{-1}$ and constants $M, \varepsilon > 0$ such that*

*(H1) For all $A \in \Lambda'$ and $n \geq 0$, $\mathbb{P}(\tau^A > n; X_1^A, \ldots, X_n^A \notin \Lambda') \leq M\rho^n$;*

*(H2) For all $A \in \Lambda'$ and $n \geq 0$, $\mathbb{P}(\tau^A > n) \leq M\mathbb{P}(\tau^{A'} > n)$;*

*(H3) For all $A \in \Lambda'$, $\mathbb{P}(X_n^A = A'$ for some $n \leq M) \geq \varepsilon$.*

*Then $P$ is $R$-positive and its left eigenvector is summable.*

To see the hypothesis (H1)-(H3) hold for the subcritical branching process, we use the following lemma.

**Lemma 3.15.** *Let $(X_t)_t$ be a subcritical branching process as in Proposition 3.1. Then, for every $\rho > 0$ there are positive constants $K$ and $M$ such that*

$$\mathbb{P}(X_1^j > K, \ldots, X_n^j > K) \leq M\rho^n$$

*Proof.* As there is $s > 0$ such that $\mathbb{E}e^{sZ} < +\infty$ then, by [NSS04, Theorem 2.1], $\mathbb{E}e^{sX_1} < +\infty$. Therefore, for all $q > 0$, $\mathbb{E}[(X_1^1)^q] < +\infty$. With this observation, the proof of the lemma is similar to the proof of [AEGR15, Proposition 3.2].

Given $\rho > 0$, there is $q$ such that $[\mathbb{E}(X_1^1)]^q < \rho$ (because $(X_k)_{k \in \mathbb{N}}$ is subcritical). By the branching property, $X_t^j$ have the same distribution as the sum of $j$ independent copies of $X_t^1$, each one of them we denote by $\tilde{X}^{(i)}$. Using the Law of Large Numbers of $L^q$,

$$\frac{X_1^j}{j} = \frac{1}{j}\sum_{i=1}^{j}\tilde{X}_1^{(i)} \xrightarrow{L^q} \mathbb{E}[X_1^1].$$

Therefore, there are positive $C$ and $K$ such that

$$\mathbb{E}\left[\frac{(X_1^1)^q}{j^q}\right] \leq C$$

for $j \leq K$ and

$$\mathbb{E}\left[\frac{(X_1^1)^q}{j^q}\right] \leq \rho$$

if $j > K$. Then, for any $1 \leq j \leq K$,

$$\mathbb{P}(X_1^j > K, X_2^j > K, \ldots, X_n^j > K) \leq \frac{1}{K^q}\mathbb{E}\left[(X_n^j)^q \mathbb{1}_{\{X_1^j > K, X_2^j > K, \ldots, X_n^j > K\}}\right]$$

$$= \frac{j^q}{K^q}\mathbb{E}\left[\frac{(X_1^j)^q}{(X_0^j)^q}\cdots\frac{(X_{n-1}^j)^q}{(X_{n-2}^j)^q}\frac{(X_n^j)}{(X_{n-1}^j)^q}\mathbb{1}_{\{X_1^j > K, X_2^j > K, \ldots, X_n^j > K\}}\right]$$

$$= \frac{j^q}{K^q}\mathbb{E}\left\{\mathbb{E}\left[\frac{(X_1^j)^q}{(X_0^j)^q}\cdots\frac{(X_{n-1}^j)^q}{(X_{n-2}^j)^q}\frac{(X_n^j)}{(X_{n-1}^j)^q}\mathbb{1}_{\{X_1^j > K, \ldots, X_n^j > K\}}\ \bigg|\ X_1, \ldots, X_{n-1}\right]\right\}$$

$$\leq \mathbb{E}\left[\frac{(X_1^j)^q}{(X_0^j)^q}\cdots\frac{(X_{n-1}^j)^q}{(X_{n-2}^j)^q}\mathbb{1}_{\{X_1^j > K, \ldots, X_{n-1}^j > K\}}\right]\cdot\sup_{\ell > K}\mathbb{E}\left[\frac{(X_n^j)^q}{(X_{n-1}^j)^q}\ \bigg|\ X_{n-1}^j = \ell\right]$$

$$\leq \rho \cdot \mathbb{E}\left[\frac{(X_1^j)^q}{(X_0^j)^q}\cdots\frac{(X_{n-1}^j)^q}{(X_{n-2}^j)^q}\mathbb{1}_{\{X_1^j>K, X_2^j>K,\ldots,X_{n-1}^j>K\}}\right].$$

By induction, we obtain

$$\mathbb{P}(X_1^j > K, X_2^j > K, \ldots, X_n^j > K) \leq \rho^{n-1}\cdot\mathbb{E}\left[\frac{(X_1^j)^q}{(X_0^j)^q}\mathbb{1}_{\{X_1^j>K\}}\right]$$

$$\leq \rho^{n-1}\cdot\sup_{\ell\leq K}\mathbb{E}\left[\frac{(X_1^\ell)^q}{(X_0^\ell)^q}\right]$$

$$\leq C\rho^{n-1}.$$

The lemma follows taking $M = \frac{C}{\rho}$. $\qquad\qquad\square$

*Proof of Proposition 3.1.* Define $(X_k)_k := \pi(\xi_k)_k$. Note that $(X_k)_k$ is a discrete chain with decay rate $R := e^\alpha$. Let $K_1 := \min\{k \geq 1 : \mathbb{P}(Z > k) > 0\}$. Fix some $\rho < R^{-1}$. Let $M$ and $K$ given by Lemma 3.15, and define $K_2 = \max\{K, K_1\}$. Hypothesis (H1) is satisfied if we take $\Lambda' = \{1, 2, \ldots, K_2\}$.

We take $A' = K_1$, and using the branching property, we have that, for $1 \leq j \leq K_2$, $\mathbb{P}(\tau^j > n) \leq \lceil\frac{j}{K_1}\rceil\mathbb{P}(\tau^{K_1} > n)$, and thus (H2) is satisfied.

As $\mathbb{P}(X_1^j = K_1) \geq (1 - e^{-1})^j e^{-K_1}[\mathbb{P}(Z = 0)]^{j-1}\mathbb{P}(Z = K_1) \geq (1 - e^{-1})^{K_2}e^{-K_1}[\mathbb{P}(Z = 0)^{K_2-1}]\mathbb{P}(Z = K_1) > 0$, (H3) is also satisfied. We conclude that $(X_k)_k$ is $R$-positive with summable left-eigenvector.

Since $\pi^{-1}\{1\} = \{\mathbf{1}\}$, $(\xi_k)_k$ is $R$-positive with left eigenvector $\nu$, therefore the continuous-time chain $(\xi_t)_t$ is $\alpha$-positive with the same eigenvectors. Summability of $\nu$ follows directly from summability of the eigenvector of the projection $(X_k)_k$

Thus, $(\xi_k)_k$ is $R$-positive with summable left eigenvector $\nu$, therefore the continuous-time chain $(\xi_t)_t$ is $\alpha$-positive with the same eigenvectors. Existence of the Yaglom limit then follows from the same arguments of [AEGR15, Theorem 3.1], that uses classical tools [SVJ66, VJ69].

Let $\nu$ and $h$ be positive vectors such that

$$\nu P_t = e^{-\alpha t}\nu; \qquad\qquad P_t h = e^{-\alpha t}h \qquad\qquad \nu h < +\infty,$$

and suppose $\nu$ is normalized to be a probability measure, i.e., $\nu\mathbb{1} = 1$, with $\mathbb{1}$ being the column vector with all entries equal to 1. Since $e^{\alpha t}\nu P_t = \nu$, for every $t \geq 0$ and $i, j \in \Lambda$,

$$e^{\alpha t} P_t(i, j) = \frac{e^{\alpha t} \nu_i P_t(i, j)}{\nu_i} \leq \frac{\nu_j}{\nu_i}.$$

Since $\Sigma_j \nu_j = 1$, by dominated convergence, we can sum over the second coordinate of (3.13) to conclude that

$$\lim_{t \to +\infty} e^{\alpha t} P_t(i, \Lambda) = h_i. \tag{3.16}$$

Hence, using (3.13) and (3.16)

$$\lim_{t \to +\infty} \mathbb{P}^i(\xi_t = j \mid \xi_t \neq \emptyset) = \lim_{t \to +\infty} \frac{P_t(i, j)}{P_t(i, \Lambda)} = \nu_i.$$

$\square$

Another way to prove existence of the Yaglom limit for the subcritical branching process is to use the fact that a discrete-time branching process with offspring distribution $\tilde{Z}$ is $R$-positive if, and only if, $\mathbb{E}[\tilde{Z} \log \tilde{Z}] < +\infty$. Moreover, in this case the right eigenvector is given by $\tilde{h}_j = j$ [SVJ66, Section 5]. Summability of the left eigenvector then follows readily from $R$-positiveness.

## 3.5 The $Q$-process

Now we prove Proposition 3.5. We make use of the following lemma, which is a version of Theorem 9 of [MV12] for countable state spaces.

**Lemma 3.17.** *Let $(\xi_t)_t$ be a Markov chain on a countable state space $\Lambda_0 = \Lambda \cup \{\emptyset\}$. Suppose that $(X_t)_t$ is absorbed at $\emptyset$ a.s. and that its sub-Markovian kernel is $\alpha$-positive with summable left eigenvector. Then, there exists a Markov process $Y = (Y_t)_t$ in $\Lambda$, whose finite dimensional distributions are given by*

$$\mathbb{P}(Y_{s_1} = i_1, \dots Y_{s_n} = i_n) = \lim_{t \to +\infty} \mathbb{P}(\xi_{s_1} = i_1, \dots, \xi_{s_n} = i_n \mid \xi_t \neq \emptyset). \tag{3.18}$$

*Moreover, $(Y_t)_t$ is time homogeneous and positive recurrent.*

*Remark 3.19.* The process $Y$ is called the *Q-process* of the process $(\xi_t)_t$.

*Proof of Lemma 3.17.* By $\alpha$-positiveness of $(\xi_t)_t$, there are positive $\nu$ and $h$ such that (3.12) holds and $\nu h = 1$. We are also assuming $\Sigma_j \nu_j = 1$.

First, we prove that the process $Y$ is well-defined, i.e., that the limit (3.18) exists. Define $\mathcal{F}_s$ to be the natural filtration of $(\xi_t)_t$. Using the Markov property of the process $(\xi_t)_t$,

$$\mathbb{P}^{i_0}(\xi_{s_1} = i_1, \ldots, \xi_{s_n} = i_n, \xi_t \neq \emptyset)$$
$$= \mathbb{E}^{i_0}(\mathbb{1}_{\{\xi_{s_1} = i_1, \ldots, \xi_{s_n} = i_n\}} \mathbb{E}^{i_0}(\mathbb{1}_{\{\xi_t \neq \emptyset\}} | \mathcal{F}_{s_n}))$$
$$= \mathbb{E}^{i_0}(\mathbb{1}_{\{\xi_{s_1} = i_1, \ldots, \xi_{s_n} = i_n\}} \mathbb{E}^{i_n}(\mathbb{1}_{\{\xi_{t-s_n} \neq \emptyset\}}))$$
$$= \mathbb{P}^{i_0}(\xi_{s_1} = i_1, \ldots, \xi_{s_n} = i_n) \mathbb{P}^{i_n}(\xi_{t-s_n} \neq \emptyset). \tag{3.20}$$

Using (3.16) and (3.20),

$$\lim_{t \to +\infty} \mathbb{P}^{i_0}(\xi_{s_1} = i_1, \ldots \xi_{s_n} = i_n \mid \xi_t \neq \emptyset)$$
$$= \lim_{t \to +\infty} \frac{\mathbb{P}^{i_0}(\xi_{s_1} = i_1, \ldots, \xi_{s_n} = i_n) \mathbb{P}^{i_n}(\xi_{t-s_n} \neq \emptyset) e^{\alpha(t-s_n)} e^{\alpha s_n}}{\mathbb{P}^{i_0}(\xi_t \neq \emptyset) e^{\alpha t}}$$
$$= \mathbb{P}^{i_0}(\xi_{s_1} = i_1, \ldots \xi_{s_n} = i_n) \frac{h_{i_n}}{h_{i_0}} e^{\alpha s_n}, \tag{3.21}$$

proving that the limit (3.18) exists.

Furthermore, the process $(Y_t)_t$ with finite dimensional distributions given by (3.18) is a Markov process. Indeed,

$$\mathbb{P}^{i_0}(Y_{s_1} = i_1, \ldots, \xi_{s_n} = i_n, Y_t = j)$$
$$= e^{\alpha t} \frac{h_j}{h_{i_0}} \mathbb{P}^{i_0}(\xi_{s_1} = i_1, \ldots, \xi_{s_n} = i_n, \xi_t = j)$$
$$= e^{\alpha(t-s_n)} e^{\alpha s_n} \frac{h_j}{h_{i_n}} \frac{h_{i_n}}{h_{i_0}} \mathbb{P}^{i_0}(\xi_{s_1} = i_1, \ldots, \xi_{s_n} = i_n) \mathbb{P}^{i_n}(\xi_{t-s_n} = j)$$
$$= \mathbb{P}^{i_0}(Y_{s_1} = i_1, \ldots, Y_{s_n} = i_n) \mathbb{P}^{i_n}(Y_{t-s_n} = j),$$

where in the first and last equalities we used (3.21) and in the second one we used the Markov property of $(\xi_t)_t$.

Therefore,

$$\mathbb{P}^{i_0}(Y_t = j \mid Y_{s_1} = i_1, \ldots, Y_{s_n} = i_n) = \mathbb{P}^{i_n}(Y_{t-s_n} = j),$$

i.e., $Y$ has the Markov property and is time homogeneous.

By (3.21) and $\alpha$-positiveness of $(\xi_t)_t$, we have that $\mathbb{P}^{i_0}(Y_t = j) \nrightarrow 0$ and thus $Y$ is positive recurrent. $\qquad\square$

We prove now Proposition 3.5. The proof consists in approximating the process $(X_t)_t$ by the process $(Y_t)_t$ using finite dimensional sets. As the process $(Y_t)_t$ is recurrent, with probability 1 it has many branching events. As the finite dimensional distributions of $X_t$ converge to those of the $Q$-process, the result follows.

*Proof of Proposition 3.5.* Denote by $Y := (Y_t)_t$ the $Q$-process of $(X_t)_t$ defined in Lemma 3.17. Let $T > 0$, its value will be determined later. Take a sequence $(\mathcal{P}_n)_n$ of partitions of $[0, T]$ such that $\mathcal{P}_{n+1} \subsetneq \mathcal{P}_n$ and $\lim_n \|\mathcal{P}_n\| = 0$. For $T > 0$, $n \in \mathbb{N}$ and $i \in \Sigma$, define the following events:

$$J_T(X) = \{\text{There are } 0 < t_1, \dots t_k < T : \ X_{t_1} = 1, X_{t_2} = 2, \dots, X_{t_{k-1}} = 1, X_{t_k} = 2\}$$

and

$$J_T^n(X) = \{\text{There are } t_1, \dots t_k \in \mathcal{P}_n : \ X_{t_1} = 1, X_{t_2} = 2, \dots, X_{t_{k-1}} = 1, X_{t_k} = 2\}.$$

We define $J_T(Y)$ and $J_T^n(Y)$ analogously.

Take $\varepsilon > 0$ arbitrary. By positive recurrence of $Y$, there is $T > 0$ such that

$$\mathbb{P}^1(J_T(Y)) \geq 1 - \varepsilon/3. \tag{3.22}$$

For the chosen $T$, there is $n_0$ such that

$$|\mathbb{P}^1((J_T^{n_0}(Y))^c) - \mathbb{P}^1((J_T(Y))^c)| \leq \varepsilon/3 \tag{3.23}$$

Finally, for all $t$ big enough, since $J_T$ and $J_T^{n_0}$ are finite dimensional cylinders, by the definition of $(Y_t)_t$,

$$|\mathbb{P}^1(J_T^{n_0}(Y))^c) - \mathbb{P}^1((J_T^{n_0}(X))^c \mid X_t \neq \emptyset)| \leq \varepsilon/3. \tag{3.24}$$

Therefore, (3.22, (3.23) and (3.24) together imply, for $t$ big enough,

$$\mathbb{P}^1((J_T(X))^c \mid X_t \neq \emptyset) \leq \mathbb{P}^1((J_T^{n_0}(X))^c \mid X_t \neq \emptyset)$$

$$\leq |\mathbb{P}^1((J_T^{n_0}(Y))^c) - \mathbb{P}^1((J_T^{n_0}(X))^c \mid X_t \neq \emptyset)|$$
$$+ |\mathbb{P}^1((J_T^{n_0}(Y))^c) - \mathbb{P}^1((J_T(Y))^c)| + \mathbb{P}^1((J_T(Y))^c)$$
$$\leq \varepsilon.$$

As $\varepsilon$ is arbitrary, this finishes the proof. $\qquad\square$